\documentclass[11pt,a4paper]{amsart}
\usepackage{amsmath, graphicx, subfigure, epsfig,amssymb,cite}
\usepackage{xcolor,mathtools}
\usepackage[normalem]{ulem}
\usepackage{amsthm,bm,mathtools}
\usepackage[a4paper,margin=1in]{geometry}

\numberwithin{equation}{section}

\usepackage{hyperref}
\hypersetup{
    colorlinks=true,
    linkcolor=blue,
    filecolor=magenta,      
    urlcolor=cyan,
}

\newcommand{\tr}{\operatorname{tr}}
\newcommand{\BP}{\mathbb P}
\newcommand{\cdn}{C_{d,n}}  

\newcommand{\E}{\mathbb{E}}
\newcommand{\Cov}{\operatorname{Cov}}

\newcommand{\Var}{\operatorname{Var}}
\newcommand{\im}{\operatorname{Im}}
\newcommand{\re}{\operatorname{Re}}
\newcommand{\e}{\epsilon}
\newcommand{\ga}{\gamma}
\newcommand{\de}{\delta}
\newcommand{\br}{\mathbb{R}}
\newcommand{\BC}{\mathbb{C}}
\newcommand{\BN}{\mathbb{N}}

\newcommand{\pa}{\partial}
\newcommand{\bt}{\beta}
\newcommand{\al}{\alpha}
\newcommand{\la}{\lambda}
\newcommand{\om}{\omega}

\newcommand{\ioi}{\int_0^{\infty}}
\newcommand{\les}{\lesssim}
\newcommand{\be}{\begin{equation}}
\newcommand{\ee}{\end{equation}}
\def\bs#1\es{
    \begin{equation}\begin{split}
    #1
    \end{split}\end{equation}
}
\newcommand{\dd}{\mathrm{d}}
\newcommand{\op}{\mathrm{op}}
\newcommand{\est}{\mathrm{e.s.t.}}

\newcommand{\CN}{\mathcal{N}}
\newcommand{\us}{\mathcal U}

\newcommand{\Eb}{\mathbb E}

\newtheorem{theorem}{Theorem}[section]
\newtheorem{lemma}[theorem]{Lemma}
\newtheorem{corollary}[theorem]{Corollary}
\newtheorem{proposition}[theorem]{Proposition}

\theoremstyle{definition}
\newtheorem{assumption}[theorem]{Assumption}

\newtheorem{remark}[theorem]{Remark}

\begin{document}

\title[High-$d$ SPA and CLT]{Saddlepoint Approximation and Central Limit Theorem for Densities in high dimensions}
\author[A Katsevich]{Alexander Katsevich$^1$}
\thanks{$^1$School of Data, Mathematical, and Statistical Sciences, University of Central Florida, Orlando, FL 32816 (Alexander.Katsevich@ucf.edu).}

\begin{abstract}
We study the saddlepoint approximation (SPA) for densities of sums of high-dimensional random vectors, allowing the dimension $d$ to grow with the sample size $n$. For averages of i.i.d.\ random vectors in $\mathbb R^d$, we obtain an explicit non-asymptotic bound for the multiplicative SPA error whose leading term is of order $d^2/n$. This improves on the generic $d^4/n$ rate in recent work of Tang and Reid.

In the natural central region $\|x\|\lesssim(d/n)^{1/2}$, we obtain a further refinement. We prove that the leading SPA error is essentially constant throughout this region, up to $O(d^{5/2}/n^{3/2})$. This allows us, at the cost of computing a single Gaussian Monte Carlo integral, to estimate the SPA error in the central region and use the resulting value as a multiplicative correction to obtain a more accurate SPA. As further corollaries to the main SPA bound, we derive a local central limit theorem with explicit multiplicative error and, on the natural high-dimensional central region, a cubic-corrected Gaussian approximation with multiplicative error $O(d^2/n)$. We also obtain a more general SPA theorem that applies, in particular, to averages of independent but not necessarily identically distributed random vectors.

We further give readily verifiable sufficient conditions for the SPA assumptions, reducing the main high-dimensional Fourier condition to one-dimensional estimates for tilted marginals. These conditions apply, in particular, to arbitrary isotropic log-concave densities, uniformly well-conditioned affine transformations thereof and their convex mixtures, as well as to distributions obtained by smoothing with Gaussian or log-concave noise. The results are illustrated by multivariate generalized Laplace distributions and by a likelihood problem for averaged images with Gaussian-mixture noise.
\end{abstract}

\maketitle

\section{Introduction}\label{sec:intro}

The saddlepoint approximation method (SPA), first proposed by Daniels in \cite{Daniels1954}, is a powerful method with numerous applications in probability and statistics \cite{jensen1995, butler2007,broda2012}. SPA provides accurate pointwise approximations to probability density functions (pdfs) and cumulative distribution functions (cdfs), based on the moment generating function (mgf). Initially proposed for scalar random variables, SPA has since been extended to vector-valued cases \cite{Kolassa2003,Kolassa2010}, primarily in fixed dimension.

The behavior of saddlepoint approximations in growing dimension has received much less attention. Jensen \cite{jensen2021} gave examples showing that the accuracy familiar from fixed-dimensional theory need not persist in high-dimensional inference. A significant recent advance is the work of Tang and Reid \cite{tang2025}, who extended Laplace and saddlepoint approximations to settings in which the dimension $d$ may grow with the sample size $n$. Under their most general assumptions, the relative SPA error is shown to be asymptotically of order $O(d^4/n)$. 

In this paper we further extend SPA analysis and obtain instead a non-asymptotic multiplicative error bound whose leading term is of order $O(d^2/n)$. Thus our result applies under the substantially weaker dimensional restriction $d^2/n\to0$. In Section~\ref{sec:VG} we consider a multivariate generalized Laplace family with uniformly bounded means and uniformly well-conditioned covariance matrices. For this family, the estimates in \cite{tang2025} give the $O(d^4/n)$ rate, while our Theorem~\ref{thm:mainI} gives $O(d^2/n)$. 

A key step in \cite{tang2025} and in our paper is an asymptotic calculation
of a high-dimensional Laplace-type integral arising after Fourier inversion
through the multivariate saddlepoint. In our notation, the relevant local
integral has the form
\be
\int_{\us}\exp\{-nG(\tau+it)\}\,\dd t,
\ee
where $\us$ is a neighborhood of zero in $\br^d$ whose size goes to zero as $d,n\to\infty$. The difference between the value of this integral, suitably normalized, and one gives the multiplicative SPA error.

Very recently, similar Laplace-type integrals have been studied by Anya Katsevich \cite{katsevich2025a}. She showed that complete asymptotic expansions can be obtained under a variety of conditions. In the case closest to ours, the leading error is also of order $O(d^2/n)$. Closely related work on the Laplace approximation of probability measures of the form $\dd\pi\propto e^{-nV(x)}\dd x$ in growing dimensions is contained in \cite{kasprzak2025,katsevich2024b}. High-dimensional posterior normality and Bernstein--von Mises-type results are also closely related; see, e.g., \cite{belloni2014,katsevich2025c} and the references therein. Such results give global Gaussian approximations of the posterior in probability metrics such as total variation, whereas the SPA studied here gives a pointwise approximation with a multiplicative error bound for the density on a prescribed region.

The approach in our paper is a blend of those in \cite{tang2025} and \cite{katsevich2025a}. Similarly to \cite{tang2025}, we use the inverse Laplace/Fourier transform and integrate across the multivariate saddlepoint; see also \cite{Kolassa2003}. Similarly to \cite{katsevich2025a}, we use a concentration inequality. This is the step that allows us to reduce the dimension dependence to $d^2$. There are two main differences from \cite{katsevich2025a}. First, the exponent in our problem is complex-valued, whereas in \cite{katsevich2025a} it is real-valued, and the argument there uses this fact in an essential way. We instead bound the oscillatory part using a Poincar\'e/spectral-gap type inequality \cite[eq.~(4.9.4)]{bakry2014}. Second, to obtain a sharper estimate, we employ thermodynamic integration, which expresses one of the key integrals in terms of a variance that can be estimated using the same inequality. The result is an explicit non-asymptotic pointwise bound for the multiplicative SPA error.

Our \emph{first main result}, Theorem~\ref{thm:mainI}, gives this bound for the density of the average of $n$ i.i.d.\ random vectors $X_i\in\br^d$. Our \emph{second main result} is an additional refinement obtained in Subsection~\ref{ssec:MC}. In centered coordinates, let $I(a)$ denote the multiplicative correction to the ordinary SPA and set $\Sigma=\Var(X)$. We show that, for every fixed $L>0$,
\be
I(a)-1=J(0)-1+O_L\big(d^{5/2}/n^{3/2}\big)+\est,
\quad
\|\Sigma^{-1/2}a\|\le L(d/n)^{1/2},
\ee
where $\est$ stands for exponentially small terms, and $J(0)$ is a single Gaussian integral evaluated at the central saddlepoint. Thus the leading SPA error is essentially constant throughout the natural high-dimensional central region. In particular, multiplying the SPA by the correction factor $J(0)$ improves the relative error from order $d^2/n$ to order $d^{5/2}/n^{3/2}$ there, a gain of a factor $(d/n)^{1/2}$; see Theorem~\ref{thm:improved SPA}. Moreover, $J(0)$ is an expectation with respect to the standard Gaussian measure and can be estimated by Monte Carlo only once for each fixed model, $d$, and $n$. Subsection~\ref{ssec:MC} provides an explicit probabilistic bound for the resulting approximation; see \eqref{MC final bnd}.

We next obtain a local central limit theorem (CLT). After affine normalization, Corollary~\ref{cor:CLT} gives an explicit multiplicative error bound. On the above natural high-dimensional central region, our \emph{third main result}, Corollary~\ref{cor:central density expansion}, further gives a cubic-corrected Gaussian approximation with multiplicative error $O(d^2/n)$. 

The correction appearing in the exponent is precisely the first multivariate Edgeworth correction; see, e.g., \cite[Chapters~2 and 4]{bhattacharya2010}. In the standard first-order Edgeworth approximation, however, this term appears additively in the polynomial prefactor multiplying the Gaussian density. Our result instead places the same correction in the exponent and can therefore be viewed as an exponentiated first-order Edgeworth approximation.

Moving this correction to the exponent is important in high dimension. A closely related phenomenon occurs for high-dimensional Laplace expansions. For the usual additive
expansion, the condition $d^2/n\to0$ is in general sharp and cannot simply be relaxed by obtaining better remainder estimates; see \cite[Example 2.19 and Proposition 2.20]{katsevich2025a}. By contrast, after passing to an expansion of the logarithm, and hence to an exponential rather than additive representation of the correction, asymptotic
expansions can remain valid in regimes arbitrarily close to the concentration threshold $d/n\to0$ \cite{katsevich2026a}. The reason is that exponentiation retains all powers generated by the leading correction, whereas an additive first-order expansion truncates after its linear contribution.

The same mechanism is relevant here. On the natural high-dimensional central region, the first Edgeworth correction is generally of order $(d^3/n)^{1/2}$ and therefore need not be small when $d^2/n\to0$. Retaining this term in the exponent nevertheless leaves a multiplicative remainder of order $O(d^2/n)$ in Corollary~\ref{cor:central density expansion}. Thus our result remains accurate in regimes in which the corresponding additive first-order Edgeworth correction is no longer perturbative. As a result we show that multiplicative convergence to the \emph{uncorrected} Gaussian density generically requires the stronger condition $d^3/n\to0$; the condition $d^2/n\to0$ is not sufficient; see \eqref{central density expansion}.

CLTs for vectors of fixed dimension are classical \cite[Sections~27 and~29]{billingsley2012}, but the case in which $d\to\infty$ together with $n$ is much less understood; see
\cite{chakraborty2025} for a recent review. An important early distributional result of Portnoy \cite{portnoy1986} identifies $d^2/n\to0$ as a fundamental scaling for multivariate central limit behavior in growing dimension.

Recent work of Bobkov and G\"otze \cite{bobkov2025} gives multidimensional Berry--Esseen bounds for the uniform local error $\sup_x |p_n(x)-\ga_d(x)|$ under substantially weaker moment assumptions than those imposed here. Their results provide global absolute bounds, while ours give pointwise multiplicative bounds. The latter can therefore capture improved accuracy in a prescribed region, including regions where the density itself is smaller than a global absolute error bound. 


Pointwise Edgeworth expansions for densities in growing dimension have also been considered recently. Withers \cite{withers2026} derives classical additive Edgeworth expansions for a broad class of statistics of dimension $d$, organized in powers of  $d^3/\sqrt n$. 

There is an extensive related literature on high-dimensional Gaussian approximation in metrics of the form $\sup_{D\in\mathcal D}|\pi(D)-\gamma(D)|$, where $\pi$ is the target distribution, $\ga$ is its Gaussian approximation, and $\mathcal D$ is a family of sets. Classical Berry--Esseen bounds for convex sets include \cite{gotze1991,bentkus2003,bentkus2005}, while sharper results are available for special geometric classes; see, for example, \cite{bentkus1996}. More recently, Gaussian approximation for hyperrectangles and related classes has been developed in \cite{chernozhukov2017} and subsequent work. Higher-order approximations for probabilities of sets in growing dimension have also been studied; see, for example, Anderson, Hall and Titterington \cite{anderson1998} for asymptotic Edgeworth expansions and Zhilova \cite{zhilova2022} for non-asymptotic Edgeworth-type expansions for Euclidean balls and half-spaces.

Many such results are obtained under substantially weaker moment assumptions and can allow the dimension to grow much faster than is possible for pointwise density approximation. This broader range of applicability reflects the fact that approximation of probabilities of sets and pointwise density approximation are fundamentally different problems, both in the form of the approximation and in the conditions required for its validity. In particular, approximation of probabilities of sets does not determine the density at an observed point. One important consequence is that it does not provide the likelihood needed for likelihood-based inference.

We also give a more general SPA theorem, Theorem~\ref{thm:mainII}, which does not require the random vector under consideration to be an average of identically distributed summands. This includes, in particular, averages of independent but non-identically distributed random vectors.

Classical SPA theory applies to broad classes of distributions in fixed dimension. In particular, uniform saddlepoint approximations for log-concave (LC) and related densities were developed by Jensen \cite{jensen1991,jensen1995}; the classical theory also includes multivariate, conditional, and non-identically distributed settings. In this paper we establish conditions under which some of these results generalize to increasing dimension $d=d_n\to\infty$.

An important feature of our results is that the assumptions are verifiable from relatively simple properties of the underlying distribution. Our \emph{fourth main result}, Lemma~\ref{lem:easy conds}, gives a sufficient criterion in terms of uniform exponential-moment bounds, covariance control, and integral H\"older regularity of one-dimensional tilted marginals. In this way the $d$-dimensional Fourier condition is reduced to 1D estimates. Thus, once the 1D tilted marginals are established to possess the required property, these multidimensional estimates need not be checked separately. Using this approach, we verify the assumptions for important and broad classes of distributions; see item (2) below.

We illustrate the theory in three complementary settings. 

\begin{enumerate}
\item Section~\ref{sec:VG} treats the multivariate generalized Laplace distribution, a non-Gaussian and generally nonseparable family for which all assumptions can be checked explicitly. It also provides the direct comparison with the dimension dependence in \cite{tang2025} discussed above. 
\item Section~\ref{sec:logconcave} applies our sufficient criterion to several broad classes of distributions. We first consider isotropic LC densities, then uniformly well-conditioned affine transformations and convex mixtures of LC densities, and then Gaussian-smoothed distributions $X=Y+Z$. Here the law of $Y$ may be supported on a lower-dimensional set, be singular, or discrete, while a uniformly nondegenerate Gaussian component $Z$ ensures the required regularity. More generally, $Z$ can be any random vector satisfying the conditions of Remark~\ref{rem:general smoothing}. 
\item Finally, Section~\ref{sec:image application} considers likelihood-based
inference from averaged images under a non-Gaussian, nonseparable
two-component Gaussian-mixture noise model.  SPAs have previously been used
to construct approximate likelihoods in fixed-dimensional and applied
multivariate problems; see, e.g., \cite{goodman2022,zhang2019_spa}.  Here we
study the same use in the growing-dimensional setting.  The exact likelihood
of the averaged image is available as a finite Gaussian mixture, allowing
direct comparison of the conventional SPA with its Gaussian and cubic
central expansions.  This illustrates a use of pointwise density
approximation that is not accessible from approximations to probabilities
of sets alone.  The experiment also tests the non-asymptotic SPA bound and
the Monte Carlo refinement of Subsection~\ref{ssec:MC}: a single Monte Carlo
estimate of the SPA error at the central saddlepoint remains representative
of the true error throughout the displayed central region.
\end{enumerate}

The paper is organized as follows. Section~\ref{sec:setting} introduces the notation and problem setting, and Section~\ref{sec:main res} proves the main SPA theorem. Section~\ref{sec:cors} contains the central Monte Carlo refinement, the local CLT, the cubic central density expansion, and the more general SPA theorem. Section~\ref{sec:VG} treats the multivariate generalized Laplace family and compares the resulting dimension dependence with \cite{tang2025}. Section~\ref{sec:sufficient criterion} develops a sufficient criterion for the assumptions, which is easier to verify. Section~\ref{sec:logconcave} applies the criterion to several broad classes of distributions. Section~\ref{sec:image application} presents the averaged-image likelihood experiment described above. Auxiliary calculations are collected in the appendices.

\section{Notation. Problem setting}\label{sec:setting}

Let $\Vert u\Vert$ denote the standard Euclidean norm in $\br^d$. 
The operator norm of a $k$-th order tensor is defined as
\be\label{T norm}
\Vert T\Vert_{\op}:=\sup_{\Vert u_1\Vert=\dots=\Vert u_k\Vert=1} \big|T[u_1,\dots, u_k]\big|,
\ee
where
\be
T[u_1,\dots, u_k]=\sum_{i_1,\dots,i_k=1}^d T_{i_1\dots i_k} (u_1)_{i_1}(u_2)_{i_2}\dots (u_k)_{i_k}.
\ee
If $T$ is symmetric, we can set $u_1=\dots=u_k$ in \eqref{T norm} \cite{zhang2012}. 
For a function $f\in C^k(\br^d)$, the $k$-th derivative tensor is
\be
(\nabla^kf(x))_{i_1\dots i_k} = \pa_{x_{i_1}}\dots\pa_{x_{i_k}}f(x). 
\ee
For matrices, we will write $\|\cdot\|$, but the operator/spectral norm is understood.

Let $X$ and $X_i$, $i\ge 1$, be independent, identically distributed (i.i.d.) random vectors in $\br^d$. 
Let the mgf and cgf of $X$ be denoted
\be\label{GFs}
\phi(\tau)=\Eb e^{\langle\tau,X\rangle},\quad \varphi(\tau)=\log\phi(\tau),
\ee
respectively. We analytically continue $\varphi(\tau)$ to $\BC^d$ in a neighborhood of the origin.
Consider a sequence of dimensions $d=d_n$, $d\to\infty$ as $n\to\infty$. Let $S^{d-1}$ be the unit sphere in $\br^d$, $I_d$ be the $d\times d$ identity matrix, and denote $\e=d/n$.

Throughout the paper we use the convention that the inequality sign $\les$ absorbs various absolute multiplicative constants which are independent of $X$, $d$, and $n$. Similarly, we write $A\les B$ for positive semidefinite matrices if $\langle Ax,x\rangle\le c \langle Bx,x\rangle$ for some $c>0$ and all $x\in\br^d$. If $A\les B$ and $B\les A$, then we write $A\asymp B$. 

Suppose $\phi(t)$ satisfies the following assumption.
\begin{assumption}\label{ass:mainI}
For all $n$ sufficiently large,
\begin{enumerate}
\item\label{exist} One has 
\be\label{Vd orig}
\sup_{\Theta\in S^{d-1}}\sup_{\tau\in V_d}\Eb \big[\langle\Theta,X\rangle^4 e^{\langle\tau,X\rangle}\big]<\infty 
\ee
for a convex neighborhood of the origin $0\in V_d\subset\br^d$.
\item\label{posdef} $H(\tau):=\nabla^2\varphi(\tau)>0$ for all $\tau\in V_d$ and $\sup_{\tau\in V_d}\Vert H(\tau)^{-1}\Vert<\infty$.
\item\label{global} One has
\bs\label{ass_1}
\bigg|\frac{\phi(\tau+i H^{-1/2} t)}{\phi(\tau)}\bigg|\le & \max\Big(\exp\big\{-\e^{1/2}\Vert  t\Vert\big\},(1+\Vert t\Vert)^{-\frac{\kappa}{d^{1/2}}}\Big),\\ 
\tau\in & V_d,\ \Vert t\Vert \ge 2.5\e^{1/2},
\es
for some $\kappa>0$. Here $H=H(\tau)$.
\item\label{phase} There exists $\de>0$ such that 
\bs
|\text{Arg}\,\phi(\tau+iH^{-1/2} t)|\le & \pi-\de,\quad |\phi(\tau+iH^{-1/2} t)|\ge\de,\\
\tau\in & V_d,\quad \Vert  t\Vert\le 2.5\e^{1/2}.
\es
\end{enumerate}
\end{assumption}
Here, $\text{Arg}\,z$ denotes the principal branch of the argument of $z\in\BC$. For notational simplicity, we drop the subscript $n$ from $d_n$ in what follows. To clarify, the constants $\kappa$ and $\de$ are independent of $d$ and $n$. In what follows we continue using the notation $H:=H(\tau)$. 

Introduce the constants
\bs\label{c3c4}
c_3:=\sup_{\tau\in V_d,\Vert t\Vert < 2.5\e^{1/2}}\Vert \nabla_t^3\im \varphi(\tau+iH^{-1/2} t)\Vert_{\op},\\
c_4:=\sup_{\tau\in V_d,\Vert t\Vert < 2.5\e^{1/2}}\Vert \nabla_t^4\re \varphi(\tau+iH^{-1/2} t)\Vert_{\op}.
\es
Assumptions~\ref{ass:mainI}(\ref{exist}, \ref{posdef}, \ref{phase}) imply that $c_3,c_4<\infty$ for each $d$, but they may grow with $d$. Also, Assumption~\ref{ass:mainI}\eqref{exist} implies that $\phi(z)$ is analytic in the polystrip $z=\tau+i t$, $\tau\in V_d$, $t\in\br^d$.

A useful stronger version of Assumption~\ref{ass:mainI}(\ref{exist}) is the
uniform bound
\be\label{uniform directional moments}
\sup_{\|\tau\|<r}\Eb e^{\langle\tau,X\rangle}\les 1
\ee
for some $r>0$ and all $d$. We have dropped here the factor $\langle\Theta,X\rangle^4$, because, 
\be\label{cosh bound}
|\langle\Theta,x\rangle|^m\le C_{h,m}(e^{h\langle\Theta,x\rangle}+e^{-h\langle\Theta,x\rangle}),\quad x\in\br^d,
\ee
for any $h,m>0$ and some $C_{h,m}$. Thus, reducing $r$ to $r'$ and taking any $0<h<r-r'$, we obtain \eqref{Vd orig} with $V_d=\{\|\tau\|\le r'\}$.

Suppose in addition that
\be\label{covariance uniform bounds}
\Var(X)\asymp I_d.
\ee
Then there exists $0<r'\le r$, independent of $d$, such that, with $V_d=\{\|\tau\|<r'\}$,
\be\label{spectral bounds}
H(\tau)\asymp I_d, \quad \tau\in V_d.
\ee
Indeed, \eqref{uniform directional moments} gives uniform bounds on the derivatives of $\phi$ of every fixed order in a fixed complex neighborhood of the origin.  After decreasing $r'$, if necessary, $\phi$ is uniformly bounded away from zero there, and hence the corresponding derivatives of $\varphi=\log\phi$ are uniformly bounded as well. In particular,
$H(0)=\Var(X)$ and $\|H(\tau)-H(0)\| \les\|\tau\|$. After decreasing $r'$ once more, if necessary, Assumption~\ref{ass:mainI}(\ref{posdef}) therefore follows from
\eqref{covariance uniform bounds} and \eqref{uniform directional moments}. The same argument shows that Assumption~\ref{ass:mainI}(\ref{phase}) holds for all $\epsilon$
sufficiently small and that the constants $c_3,c_4$ in \eqref{c3c4}
are bounded independently of $d$.

We summarize the preceding discussion in a lemma.
\begin{lemma}\label{lem:aux ass}
Let the random vector $X$ satisfy \eqref{uniform directional moments} and \eqref{covariance uniform bounds}. There exist $0<r'\le r$ and $\e_0>0$ such that, provided $\e=d/n\le\e_0$, 
\begin{enumerate}
\item parts \eqref{exist}, \eqref{posdef} and \eqref{phase} of Assumption~\ref{ass:mainI} are satisfied with $V_d:=\{\|\tau\|<r'\}$, and 
\item constants $c_3,c_4$, defined in \eqref{c3c4} using the same $V_d$, are bounded independently of $d$, $n$. 
\end{enumerate}
\end{lemma}

Introduce the key function 
\bs\label{G def}
G_a(t)=&-\varphi(\tau+i H^{-1/2}t)+\varphi(\tau)+i\langle H^{-1/2}t,a\rangle,\\ \tau=&(\nabla\varphi)^{-1}(a),\ \Vert t\Vert < 2.5\e^{1/2},\ a\in U_d:=\nabla\varphi(V_d). 
\es
This relation between $a$ and $\tau$ is assumed in what follows. Note that $(\nabla\varphi)^{-1}(a)$ is well-defined (i.e., single-valued) due to the assumption $\nabla^2\varphi(\tau)>0$, $\tau\in V_d$. By assumption and construction, 
\be\label{G t-exp}
G_a(t)=\tfrac12\Vert t\Vert^2 +O(\Vert t\Vert^3),\ t\to0.
\ee
The absence of the linear term in the Taylor expansion follows from the relation $\nabla\varphi(\tau)=a$ and the analyticity of $\varphi(z)$ in a neighborhood of $z=\tau$. Indeed, since $\varphi(z)$ is analytic, the derivative can be taken along any direction in $\BC^d$. Therefore,
\be\label{cgf exp}
\varphi(\tau+i H^{-1/2}t)=\varphi(\tau)+\nabla \varphi(\tau)i H^{-1/2}t-\tfrac12\Vert t\Vert^2+O(\Vert t\Vert^3),\ t\to0,
\ee
where $\nabla \varphi(\tau)=\nabla_\tau \varphi(\tau)$, and the assertion follows.

From \eqref{GFs} it follows immediately that $\re G_a(t)$ and $\im G_a(t)$ are even and odd in $t$, respectively.

Let us now discuss Assumption~\ref{ass:mainI}\eqref{global}. The bound on the right is the maximum of two functions. The first of these is similar to the one in \cite[Assumption 2.3]{katsevich2025a}. This function dominates on the disk $\|t\|\le \rho_n$, where $\rho_n$ solves the equation
\be\label{crossover}
\frac{\log(1+\rho_n)}{\rho_n}=\frac{d}{\kappa\sqrt n}.
\ee
Thus, $\rho_n$ grows slightly faster than linearly with $\sqrt n/d$. Taking logarithm on both sides of \eqref{ass_1} we get
\be\label{local cond v1}
\re\varphi(\tau+iH^{-1/2} t)-\varphi(\tau)\le -\e^{1/2}\Vert  t\Vert,\quad  2.5\e^{1/2}\le \Vert t\Vert\le \rho_n.
\ee
Using \eqref{cgf exp}, this means that for any fixed $c>2.5$, one should have
\be\label{local cond v2}
\e^{1/2}\Vert  t\Vert\le \tfrac12\Vert t\Vert^2+O(\|t\|^3),\quad 2.5\e^{1/2}\le \Vert t\Vert \le c\e^{1/2},
\ee
which clearly holds provided $\e$ is sufficiently small. This means that the first term in \eqref{ass_1} is consistent with the local quadratic behavior in \eqref{cgf exp}. See \cite[Remark 2.4]{katsevich2025a} for a more detailed discussion of this condition.

For $\Vert t\Vert\ge\rho_n$, the second function inside the maximum in \eqref{ass_1} dominates. It is used to avoid exponential decay at infinity. In fact, since $d\to\infty$, we see that a fairly slow decay as $t\to\infty$ is allowed when $d\gg1$.


Observe that the mgf of $\bar X_n:=(1/n)\sum_{i=1}^n X_i$ is $\phi^n(\tau/n)$. By the Laplace transform inversion formula, the pdf of $\bar X_n$ is given by
\be
\rho_n(y)=\frac1{(2\pi)^d}\int_{\br^d}\phi^n((p+iq)/n)e^{-\langle p+iq,y\rangle}\dd q,\quad y\in\br^d,
\ee
for any $p$ such that $p/n\in V_d$. Pick any $a\in U_d$, set $y=a$, select $p$ so that $\tau=p/n\in V_d$, and change variables $t=H^{1/2}q/n$. This gives
\bs
\rho_n(a)&=\frac{n^d}{(2\pi)^d\sqrt{\det H}}\int_{\br^d}\phi^n\big(\tau+iH^{-1/2}t\big)e^{-n\langle \tau+iH^{-1/2}t,a\rangle}\dd t\\
&=\frac{n^d}{(2\pi)^d\sqrt{\det H}}e^{-n\varphi^*(a)}\int_{\br^d}\bigg[\frac{\phi\big(\tau+iH^{-1/2}t\big)}{\phi(\tau)}\bigg]^n e^{-in\langle H^{-1/2}t,a\rangle}\dd t,\\
\varphi^*(a)&=\langle\tau,a\rangle-\varphi(\tau),\quad a\in U_d.
\es
Denote
\bs\label{I(a) def}
I(a)=&\cdn\int_{\br^d}\bigg[\frac{\phi\big(\tau+iH^{-1/2}t\big)}{\phi(\tau)}\bigg]^n e^{-in\langle H^{-1/2}t,a\rangle}\dd t,\quad \cdn:=\frac{n^{d/2}}{(2\pi)^{d/2}}.
\es
Then
\bs\label{rho a}
\rho_n(a)&=\frac{\cdn}{\sqrt{\det H}} e^{-n\varphi^*(a)}I(a),\quad a\in U_d.
\es
The conventional SPA step consists of replacing $I(a)$ by 1 in \eqref{rho a}.
%
%
%
%

\section{First SPA result and its proof}\label{sec:main res}

Set 
\bs\label{lots of notation}
&R=2.5,\quad 
\us:=\{t\in\br^d:\Vert t\Vert< R\e^{1/2}\},\quad
\widetilde\us:=\{s\in\br^d:\Vert s\Vert< Rd^{1/2}\},\\
&\la_*:=1-\frac{c_4R^2}2\frac dn,\quad
A_0:=\frac{c_4d(d+2)}{24n},\quad
B_0:=\frac{c_4^2d(d+4)(d+9)}{72\la_*^4n^2},\\
&D_0:=A_0+B_0.
\es

\begin{theorem}\label{thm:mainI}
Suppose Assumption~\ref{ass:mainI} is satisfied and $\la_*>0$. For all
$d^2/n$ sufficiently small, one has
\bs\label{main res I}
|I(a)-1|\le &
e^{D_0}\bigg(D_0+\frac{c_3^2}{8\la_*^3}\frac{d(d+4)}n\bigg)
+c\left[\exp(-d)+\bigg(\frac e{\kappa^2}\frac{d^2}n\bigg)^{d/2}\right],\
a\in U_d,
\es
for an absolute constant $c>0$.
\end{theorem}

Importantly, the unspecified absolute constant $c$ appears only in the exponentially small remainder. It can be computed, but we do not do it to avoid irrelevant technical details. The leading, algebraically decaying part of the bound is completely explicit and has a simple, transparent form in terms of $c_3$, $c_4$, and $\lambda_*$.  In particular, if $c_3,c_4=O(1)$ and $d^2/n\to0$, then
\be
\la_*^{-3}=1+O(d/n),\quad 
D_0=\Big(\frac{c_4}{24}+O(1/d)\Big)\frac{d^2}n,
\ee
and we get an estimate
\be\label{main res leading}
|I(a)-1|\le
\exp\Big\{\big(\tfrac{c_4}{24}+O(1/d)\big)\frac{d^2}n\Big\}\Big(\frac{3c_3^2+c_4}{24}+O(1/d)\Big)\frac{d^2}{n}+\est\,.
\ee
Here and throughout the paper, $\est$ denotes a term bounded in absolute value, up to a constant prefactor independent of $d$ and $n$, by
\be
\mathcal E_{d,n}:=e^{-d}+\left(\frac{e}{\kappa^2}\frac{d^2}{n}\right)^{d/2}.
\ee

Our proof broadly follows the approach in \cite{katsevich2025a}. The tail estimate is treated similarly, but there are two important differences in the analysis of the local integral.  First, the exponent $-n G_a(t)$ in our problem (see \eqref{Iu}) is complex-valued, whereas the exponent in \cite{katsevich2025a} is real-valued. The logarithmic Sobolev argument used there therefore does not apply directly to the oscillatory factor. We instead separate the real and imaginary parts of the exponent and bound the contribution of the phase using a Poincar\'e/Brascamp--Lieb inequality for the probability measure generated by $\exp\{-n\re G_a(t)\}$.

Second, to carefully estimate one of the integrals involving $\exp\{-n\re G_a(t)\}$, we use thermodynamic integration. This expresses the remainder in terms of a variance, which can again be bounded by the same Poincar\'e/Brascamp--Lieb inequality.

\subsection{Tail integrals}\label{ssec:tail}

Denote (cf. \eqref{I(a) def})
\be\label{Ic def}
I_{\us^c}:=\cdn\int_{\us^c}
\bigg[\frac{\phi\big(\tau+iH^{-1/2}t\big)}{\phi(\tau)}\bigg]^n
e^{-in\langle H^{-1/2}t,a\rangle}\dd t.
\ee
By \eqref{ass_1},
\bs\label{prep Ic}
|I_{\us^c}|\le & I_1+I_2,\\
I_1:=&\cdn\int_{\Vert t\Vert\ge R\e^{1/2}}
\exp\big\{-(dn)^{1/2}\Vert t\Vert\big\}\dd t,\\
I_2:=&\cdn\int_{\br^d}(1+\Vert t\Vert)^{-\mu}\dd t,\;
\mu:=\kappa\frac{n}{d^{1/2}}.
\es

Write the integral in $I_1$ in spherical coordinates and change variables
\bs\label{I1 v1}
I_1=&\frac{|S^{d-1}|}{(2\pi)^{d/2}}
\int_{R d^{1/2}}^\infty
\exp\big\{-d^{1/2}\rho\big\}\rho^{d-1}\dd\rho\\
=&\frac{|S^{d-1}|}{(2\pi)^{d/2}}
\int_{R d^{1/2}}^\infty\exp\{-\psi(\rho)\}\dd\rho,\\
\psi(\rho):=&d^{1/2}\rho-(d-1)\log\rho.
\es
Here $|S^{d-1}|$ is the area of the unit sphere in $\br^d$. Since $R>1$, the critical point
$\rho_*=(d-1)/d^{1/2}$ satisfies $\rho_*<Rd^{1/2}$. Hence
\be
\psi^\prime(\rho)\ge\psi^\prime(Rd^{1/2})
\ge d^{1/2}[1-(1/R)],\; \rho\ge Rd^{1/2},
\ee
and therefore
\be
\psi(\rho)\ge\psi(Rd^{1/2})
+d^{1/2}[1-(1/R)](\rho-Rd^{1/2}),\; \rho\ge Rd^{1/2}.
\ee
Using this inequality in \eqref{I1 v1}, integrating, and using the Stirling formula, we obtain
\bs\label{I1 bnd_1}
I_1\le&
\frac{|S^{d-1}|}{(2\pi)^{d/2}}
\frac{R^{d-1}d^{(d-1)/2}e^{-Rd}}
{d^{1/2}[1-(1/R)]}
\les 
\frac{R^de^{d/2}e^{-Rd}}{d^{1/2}(R-1)}.
\es
Since $R=2.5$ satisfies $R\ge[(1/2)+\log R]+1$,
\be\label{I1 bnd_2}
I_1\les\frac{\exp(-d)}{d^{1/2}}.
\ee

Finally, consider $I_2$. Our assumptions imply
\be
d/\mu=\frac{d^{3/2}}{\kappa n}
=\frac{1}{\kappa d^{1/2}}\frac{d^2}{n}\to0
\ee
as $d,n\to\infty$. Evaluate the integral in the definition of $I_2$ and use the Stirling formula to get
\bs
I_2=&\cdn|S^{d-1}|\ioi(1+\rho)^{-\mu}\rho^{d-1}\dd\rho\\
=&\frac{n^{d/2}}{(2\pi)^{d/2}}
\frac{2\pi^{d/2}}{\Gamma(d/2)}
\frac{\Gamma(d)\Gamma(\mu-d)}{\Gamma(\mu)}
=\frac{(2n)^{d/2}\Gamma((d+1)/2)}{\pi^{1/2}}
\frac{\Gamma(\mu-d)}{\Gamma(\mu)}.
\es
Using Stirling's formula we compute
\be
I_2\les
\frac{(edn)^{d/2}(\mu-d)^{\mu-d}}{\mu^\mu}
=\bigg(\frac{edn}{\mu^2}\bigg)^{d/2}
(1-(d/\mu))^{\mu-d}.
\ee
Since $edn/\mu^2=(e/\kappa^2)d^2/n$ and
$(1-(d/\mu))^{\mu-d}<1$,
\be
I_2\les
\bigg(\frac e{\kappa^2}\frac{d^2}{n}\bigg)^{d/2}.
\ee

\subsection{Local integrals}\label{ssec:local}

Consider the part of the integral in \eqref{I(a) def} over $\us$.
Assumption~\ref{ass:mainI}\eqref{phase} implies that the cgf
$\varphi(\tau+iH^{-1/2}t)$ is well defined for $\tau\in V_d$,
$t\in\us$. Using \eqref{G def},
\be\label{Iu}
I_{\us}=\cdn\int_{\us}e^{-nG_a(t)}\dd t.
\ee
For simplicity, drop the subscript $a$ from $G_a$ and write
$G(t)=G_r(t)+iG_i(t)$. By construction,
\be\label{Gri exp}
G_r(t)=\tfrac12\Vert t\Vert^2+O(\Vert t\Vert^4),\;
G_i(t)=O(\Vert t\Vert^3),\; t\to0.
\ee
Set $r(t):=G_r(t)-\tfrac12\Vert t\Vert^2$. Changing variables
$s=n^{1/2}t$ gives
\bs\label{local decomposition new}
I_{\us}=&I_0-I_r-iI_i,\\
I_0:=&\int_{\widetilde\us}e^{-q(s)}\ga(\dd s),\;
q(s):=nr(s/n^{1/2}),\\
I_r:=&\int_{\widetilde\us}
\big(1-\cos(nG_i(s/n^{1/2}))\big)e^{-q(s)}\ga(\dd s),\\
I_i:=&C_{d,n}\int_{\us}
\sin(nG_i(t))e^{-n G_r(t)}\dd t,
\es
where $\ga=\CN(0,I_d)$ denotes the standard Gaussian measure on $\br^d$. As was mentioned before, $G_r$ is even and $G_i$ is odd. Using that $\us$ is symmetric gives $I_i=0$.

To bound the remaining integrals, we first obtain a convexity estimate which will be used repeatedly. Since $G_r$ is even, $\nabla^2G_r(0)=I_d$, and $\nabla^3G_r(0)=0$, we get
\be\label{Hess Gr expansion}
\nabla^2G_r(t)-I_d
=\int_0^1(1-u)\nabla^4G_r(ut)[t,t,\cdot,\cdot]\dd u.
\ee
By \eqref{c3c4}, for $t\in\us$,
\be\label{Hess Gr lower}
\Vert\nabla^2G_r(t)-I_d\Vert
\le\tfrac{c_4}{2}\Vert t\Vert^2\le \tfrac{c_4R^2}2\tfrac dn,
\ee
and therefore
\be\label{Gr strong convexity}
\nabla^2G_r(t)\succeq\la_*I_d,\quad t\in\us.
\ee

We use the following result. If $\nu$ is a probability measure
on $\widetilde\us$ with density proportional to $e^{-V(s)}$ and
$\nabla^2V\succeq\la_*I_d$, then by \cite[Proposition 2.1]{bobkov2000}
\be\label{BL new}
\Var_\nu(f)\le
\tfrac1{\la_*}\int_{\widetilde\us}\Vert\nabla f(s)\Vert^2\nu(\dd s).
\ee
For the measures $\nu$ considered below, $\E_\nu s=0$. Hence, applying \eqref{BL new} first to the coordinate functions $f=s_j$ and then to $f=\|s\|^2$ gives
\be\label{moments 2 4 new}
\int_{\widetilde\us}\Vert s\Vert^2\nu(\dd s)\le\tfrac d{\la_*},\quad
\int_{\widetilde\us}\Vert s\Vert^4\nu(\dd s)
\le\tfrac{d(d+4)}{\la_*^2}.
\ee
Similarly, applying \eqref{BL new} to $\Vert s\Vert^3$ and using
Cauchy--Schwarz,
\be\label{moment 6 new}
\int_{\widetilde\us}\Vert s\Vert^6\nu(\dd s)
\le\frac{d(d+4)(d+9)}{\la_*^3}.
\ee

We next estimate $I_r$. Define the probability measure
\be\label{nu1 def}
\nu_1(\dd s):=
I_0^{-1}e^{-q(s)}1_{\widetilde\us}(s)\ga(\dd s)
\ee
and let $\upsilon(s):=\sin(nG_i(s/n^{1/2})/2)$. Here we use that $I_0$, defined in \eqref{local decomposition new}, is the normalizing constant of the measure proportional to
$e^{-q(s)}1_{\widetilde\us}(s)\ga(\dd s)$. Since $\upsilon$ is odd and $q$ in \eqref{local decomposition new} is even, we see that $\nu_1$ is even and $\Eb_{\nu_1}\upsilon=0$. Thus
\be\label{Ir variance new}
I_r=2I_0\Eb_{\nu_1}\upsilon^2
=2I_0\Var_{\nu_1}(\upsilon).
\ee
The potential of $\nu_1$, apart from an additive constant, is
\be
V_1(s)=\tfrac12\Vert s\Vert^2+q(s)
=nG_r(s/n^{1/2}).
\ee
By \eqref{Gr strong convexity},
$\nabla^2V_1(s)\succeq\la_*I_d$ on $\widetilde\us$. Therefore, by
\eqref{BL new}, with $\nu=\nu_1$:
\be\label{Ir BL new}
I_r\le
\frac{2I_0}{\la_*}
\int_{\widetilde\us}\Vert\nabla\upsilon(s)\Vert^2\nu_1(\dd s).
\ee
We have
\be
\nabla\upsilon(s)=
\frac{n^{1/2}}2
\cos\big(nG_i(s/n^{1/2})/2\big)
\nabla G_i(s/n^{1/2}).
\ee
Using \eqref{grad bnd} and \eqref{c3c4},
\be\label{grad upsilon new}
\Vert\nabla\upsilon(s)\Vert
\le\frac{c_3}{4n^{1/2}}\Vert s\Vert^2.
\ee
Combining \eqref{Ir BL new}, \eqref{grad upsilon new}, and
\eqref{moments 2 4 new},
\be\label{Ir new bnd}
I_r\le
\frac{c_3^2I_0}{8\la_*^3n}d(d+4).
\ee

It remains to estimate $I_0$. For $0\le\bt\le1$, set
\be\label{thermo Z}
Z(\bt):=\int_{\widetilde\us}e^{-\bt q(s)}\ga(\dd s),\quad
\nu_\bt(\dd s):=
Z(\bt)^{-1}e^{-\bt q(s)}1_{\widetilde\us}(s)\ga(\dd s).
\ee
Thus $Z(0)=\ga(\widetilde\us)$ and $Z(1)=I_0$. Differentiating under the integral,
\be\label{thermo derivs}
(\log Z)^\prime(\bt)=-\Eb_{\nu_\bt}q,\quad
(\log Z)^{\prime\prime}(\bt)=\Var_{\nu_\bt}(q).
\ee
Taylor's formula with integral remainder gives the identity
\be\label{thermo identity}
\log(I_0/\ga(\widetilde\us))
=-\Eb_{\nu_0}q
+\int_0^1(1-\bt)\Var_{\nu_\bt}(q)\dd\bt.
\ee

The first term is estimated directly. Since $\nabla^l r(0)=0$, $0\le l\le 3$, Taylor's formula gives
\be\label{r fourth exact}
r(t)=\frac16\int_0^1(1-u)^3
\nabla^4G_r(ut)[t,t,t,t]\dd u.
\ee
Since $\nu_0$ is the standard Gaussian measure conditioned on the centered ball $\widetilde\us$, conditioning can only decrease the fourth radial moment. Hence from \eqref{local decomposition new} and \eqref{r fourth exact}
\be\label{mean q bnd}
|\Eb_{\nu_0}q|
\le\frac{c_4}{24n}\Eb_{\nu_0}\Vert s\Vert^4
\le\frac{c_4}{24n}\Eb_{\ga}\Vert s\Vert^4
\le\frac{c_4}{24n}d(d+2)=A_0.
\ee

For the variance term in \eqref{thermo identity}, the potential of $\nu_\bt$ is
\be
V_\bt(s)=\tfrac12\Vert s\Vert^2+\bt q(s).
\ee
Since
$q(s)=nG_r(s/n^{1/2})-\tfrac12\Vert s\Vert^2$,
\be
\nabla^2V_\bt(s)
=(1-\bt)I_d+\bt\nabla^2G_r(s/n^{1/2})
\succeq\la_*I_d;
\ee
see \eqref{Gr strong convexity}. Also, Taylor's formula applied to $\nabla r$ gives
\be\label{grad r fourth}
\Vert\nabla r(t)\Vert\le
\frac{c_4}{6}\Vert t\Vert^3,
\ee
and hence
\be\label{grad q new}
\Vert\nabla q(s)\Vert\le
\frac{c_4}{6n}\Vert s\Vert^3.
\ee
Applying \eqref{BL new} to $q$ under $\nu_\bt$ and using
\eqref{moment 6 new},
\be\label{var q new}
\Var_{\nu_\bt}(q)
\le\frac{c_4^2}{36\la_*n^2}
\Eb_{\nu_\bt}\Vert s\Vert^6
\le
\frac{c_4^2}{36\la_*^4n^2}d(d+4)(d+9).
\ee
Substituting \eqref{mean q bnd} and \eqref{var q new} into
\eqref{thermo identity},
\be\label{log I0 bnd new}
|\log(I_0/\ga(\widetilde\us))|\le A_0+B_0=D_0.
\ee
Consequently, $I_0\le e^{D_0}$ and
\be\label{I0 new bnd}
|I_0-1|\le(1-\ga(\widetilde\us))+D_0e^{D_0}.
\ee

Combining \eqref{Ir new bnd} and \eqref{I0 new bnd}, and using
$I_0\le e^{D_0}$,
\be\label{local final new}
|I_{\us}-1|
\le(1-\ga(\widetilde\us))
+e^{D_0}\bigg(
D_0+\frac{c_3^2}{8\la_*^3}\frac{d(d+4)}n
\bigg).
\ee
The standard Gaussian tail bound used above gives $1-\ga(\widetilde\us)\le\exp\{-(R-1)^2d/2\}$. Since $R=2.5$, this term is absorbed by $\exp(-d)$. Combining \eqref{local final new} with the estimates for $I_1$ and $I_2$ proves Theorem~\ref{thm:mainI}.

\section{Consequences and refinements}\label{sec:cors}

\subsection{A Monte Carlo refinement in the central region}
\label{ssec:MC}

The bound in Theorem~\ref{thm:mainI} is uniform in $a\in U_d$. For numerical purposes, a sharper bound can be obtained if $a$ is restricted to a central region. We describe this observation in centered coordinates, so that $\Eb X=0$ and hence $\tau(0)= (\nabla\varphi)^{-1}(0)=0$. The noncentered case follows by translation, with the central region centered at $\Eb X$.

In this subsection we suppose that \eqref{uniform directional moments} and
\eqref{covariance uniform bounds} hold. Let $V_d=\{\|\tau\|<r'\}$, with $r'>0$ chosen as in Lemma~\ref{lem:aux ass}, and suppose also that \eqref{ass_1} holds on this $V_d$.

Set $\Sigma:=H(0)=\Var(X)$. Under the Gaussian approximation, the covariance of the sample mean is $\Sigma/n$, and the typical standardized radius is of order $\sqrt d$. Accordingly, for fixed $L>0$ we take the central region to be
\be\label{central region}
\|\Sigma^{-1/2}a\|\le L(d/n)^{1/2}.
\ee
Under the uniform spectral bounds \eqref{spectral bounds}, this is equivalent, up to changing $L$, to $\|a\|\le L(d/n)^{1/2}$. For every fixed $L>0$, these same bounds imply
that, for $d/n$ sufficiently small, the central region \eqref{central region} is contained in $U_d$. Moreover, $\|\tau(a)\|\les_L (d/n)^{1/2}$ uniformly on this region.

Restore temporarily the dependence on $a$ in the notation of Subsection~\ref{ssec:local}. For $\tau=\tau(a)$, set
\be\label{MC J}
J(a)
:=
\int_{\widetilde\us}
e^{-q_\tau(s)}\cos\theta_\tau(s)\,\gamma(\dd s),
\quad
q_\tau(s)
:=
nG_{r,a}(s/\sqrt n)-\tfrac12\|s\|^2,
\quad
\theta_\tau(s)
:=
nG_{i,a}(s/\sqrt n).
\ee
Thus $J(a)=I_0(a)-I_r(a)$. Since $I_i=0$, the tail estimates in Subsection~\ref{ssec:tail} give $I(a)=J(a)+\est$.

Fix any $L>0$. 
Since $G_{r,a}(t)-\|t\|^2/2$ is even and begins at fourth order in $t$, while $G_{i,a}(t)$ is odd and begins at third order, Taylor expansion in $\tau$ about $\tau=0$ gives
\be\label{MC Taylor}
q_\tau(s)=q_0(s)+R_q(s),
\quad
\theta_\tau(s)=\theta_0(s)+R_\theta(s),
\ee
where, uniformly for $s\in\widetilde\us$,
\be\label{MC remainder bounds}
|R_q(s)|
 \les_L
\frac{\sqrt d\,\|s\|^4}{n^{3/2}},
\quad
|R_\theta(s)|
 \les_L
\frac{\sqrt d\,\|s\|^3}{n}.
\ee
Both bounds tend to zero uniformly on $\widetilde\us$ provided that $d^2/n\to0$. The first estimate uses derivatives through order five, whereas the second requires only derivatives through order four (cf. the paragraph following \eqref{spectral bounds}).

Substituting \eqref{MC Taylor} into \eqref{MC J} gives
\be
J(a)
=
\int_{\widetilde\us}
e^{-q_0(s)}e^{-R_q(s)}
\cos\bigl(\theta_0(s)+R_\theta(s)\bigr)\ga(\dd s).
\ee
Since $d^2/n$ is sufficiently small, $q_0$ and $R_q$ are uniformly bounded on $\widetilde\us$. Hence
\be
e^{-R_q}=1-R_q+O(R_q^2),
\quad
\cos(\theta_0+R_\theta)
=
\cos\theta_0
-
R_\theta\sin\theta_0
+
O(R_\theta^2),
\ee
and therefore
\be\label{MC integrand expansion}
J(a)-J(0)
=
-\int_{\widetilde\us}e^{-q_0}R_q\cos\theta_0\dd\ga
-\int_{\widetilde\us}e^{-q_0}R_\theta\sin\theta_0\dd\ga
+
O\left(
\int_{\widetilde\us}
(R_q^2+R_\theta^2)\dd\ga\right).
\ee

The first integral is bounded directly using \eqref{MC remainder bounds}:
\be
\int_{\widetilde\us}|R_q|\dd\ga
\les_L
\frac{\sqrt d\,\Eb\|S\|^4}{n^{3/2}}
\les_L
\frac{d^{5/2}}{n^{3/2}},
\quad S\sim \CN(0,I_d).
\ee

For the second integral we use the fact that both $\theta_0$ and $R_\theta$ are odd. Moreover,
\be
\|\nabla_s\theta_0(s)\|
\les
\frac{\|s\|^2}{\sqrt n},
\quad
\|\nabla_sR_\theta(s)\|
\les_L
\frac{\sqrt d\,\|s\|^2}{n}.
\ee
Applying \eqref{BL new} to the Gaussian measure conditioned on the symmetric ball $\widetilde\us$ therefore gives
\be
\|\theta_0\|_{L^2(\ga_{\widetilde\us})}
 \les
\frac d{\sqrt n},
\quad
\|R_\theta\|_{L^2(\ga_{\widetilde\us})}
 \les_L
\frac{d^{3/2}}n.
\ee
Consequently, since $|\sin x|\le |x|$,
\be
\int_{\widetilde\us}
|R_\theta\sin\theta_0|\dd\ga
\les
\Big[
\int_{\widetilde\us}R_\theta^2\dd\ga
\int_{\widetilde\us}\theta_0^2\dd\ga
\Big]^{1/2}
\les_L
\frac{d^{5/2}}{n^{3/2}}.
\ee
Finally,
\be
\int_{\widetilde\us}R_\theta^2\dd\ga
\les_L\frac{d^3}{n^2},
\quad
\int_{\widetilde\us}R_q^2\dd\ga
\les_L\frac{d^5}{n^3},
\ee
which are smaller than $d^{5/2}/n^{3/2}$ when $d^2/n$ is sufficiently small. Thus \eqref{MC integrand expansion} yields
\be\label{MC center}
J(a)
=
J(0)
+
O_L\bigl(d^{5/2}/n^{3/2}\bigr),
\quad
\|\Sigma^{-1/2}a\|\le L(d/n)^{1/2}.
\ee
Hence
\be\label{MC center I}
I(a)
=
J(0)
+
O_L\bigl(d^{5/2}/n^{3/2}\bigr)
+
\est,
\quad
\|\Sigma^{-1/2}a\|\le L(d/n)^{1/2}.
\ee

This observation has two useful consequences. First, the ordinary saddlepoint approximation amounts to replacing $I(a)$ by $1$. Therefore,
\be\label{MC SPA error}
I(a)-1
=
J(0)-1
+
O_L\bigl(d^{5/2}/n^{3/2}\bigr)
+
\est,
\quad
\|\Sigma^{-1/2}a\|\le L(d/n)^{1/2}.
\ee
Thus the quantity $J(0)-1$ describes the leading SPA error throughout the central region with accuracy
\be
O_L\left(\frac{d^{5/2}}{n^{3/2}}\right)
=
O_L\left(\frac{d^2}{n}\sqrt{\frac dn}\right),
\ee
which is smaller by a factor of $(d/n)^{1/2}$ than the generic $O(d^2/n)$ SPA error scale.

Second, either $I(0)$ or $J(0)$ can be used as a multiplicative correction factor. Since $I(0)=J(0)+\est$, \eqref{MC center I} gives
\be\label{MC calibrated I}
I(a)
=
I(0)
+
O_L\bigl(d^{5/2}/n^{3/2}\bigr)
+
\est
\ee
uniformly for $\|\Sigma^{-1/2}a\|\le L(d/n)^{1/2}$.

Since $I(0)$, and hence also $J(0)$, is bounded away from zero for sufficiently small $d^2/n$ by Theorem~\ref{thm:mainI}, the preceding argument proves the following result.

\begin{theorem}[Improved saddlepoint approximation]
\label{thm:improved SPA}
Suppose that \eqref{uniform directional moments} and \eqref{covariance uniform bounds} hold. Let $V_d=\{\|\tau\|<r'\}$, with $r'>0$ chosen as in Lemma~\ref{lem:aux ass}, and suppose that \eqref{ass_1} holds on this $V_d$. Suppose also that $\Eb X=0$. Let $K_0$ be either $I(0)$ or $J(0)$. Then, for every fixed $L>0$, provided $d^2/n$ is sufficiently small,
\be\label{improved I}
\sup_{\|a\|\le L(d/n)^{1/2}}
\left|
\frac{I(a)}{K_0}-1
\right|
 \les_L
\frac{d^{5/2}}{n^{3/2}}
+
\est.
\ee
Consequently, if
\be\label{calibrated SPA}
\rho_n^{\rm imp}(a)
:=
\frac{C_{d,n}K_0}{\sqrt{\det H(\tau)}}
e^{-n\varphi^*(a)},
\quad
\tau=(\nabla\varphi)^{-1}(a),
\ee
then
\be\label{improved SPA}
\sup_{\|a\|\le L(d/n)^{1/2}}
\left|
\frac{\rho_n(a)}{\rho_n^{\rm imp}(a)}-1
\right|
\les_L
\frac{d^{5/2}}{n^{3/2}}
+
\est.
\ee
\end{theorem}

The representation \eqref{MC J} is especially convenient for Monte Carlo computation, since the reference measure is standard Gaussian. Let $S_1,\ldots,S_M$ be independent $\CN(0,I_d)$ random vectors and define
\be\label{MC estimator}
\widehat J_M(0)
:=
\frac1M
\sum_{k=1}^M
\chi_{\widetilde\us}(S_k)
e^{-q_0(S_k)}
\cos\theta_0(S_k).
\ee
Here $\chi_{\widetilde\us}$ is the characteristic function of the set $\widetilde\us$. Then $\Eb\widehat J_M(0)=J(0)$. Moreover, set
\be\label{MC B}
B_{d,n}
:=
\exp\left\{
\frac{c_4R^4}{24}\frac{d^2}{n}
\right\}.
\ee
By the bound on $q_0$, every summand in \eqref{MC estimator} is bounded in absolute value by $B_{d,n}$. Hoeffding's inequality \cite[Theorem 1]{mauer2009}, \cite[Theorem 2.2.6]{vershynin2018} therefore gives, with probability at least $1-\eta$ over the Monte Carlo sample,
\be\label{MC probabilistic bnd}
|\widehat J_M(0)-J(0)|
\le
B_{d,n}
\sqrt{\frac{2\log(2/\eta)}{M}}.
\ee
Combining this with \eqref{MC center I}, we obtain with the same probability
\be\label{MC final bnd}
\sup_{\|\Sigma^{-1/2}a\|\le L(d/n)^{1/2}}
|I(a)-\widehat J_M(0)|
 \les_L
\frac{d^{5/2}}{n^{3/2}}
+
B_{d,n}
\sqrt{\frac{2\log(2/\eta)}{M}}
+
\est.
\ee
When $c_4 \les1$ and $d^2/n$ is sufficiently small, $B_{d,n} \les1$. Thus the Monte Carlo error has the usual $M^{-1/2}$ rate with a dimension-independent prefactor. The quantity $\widehat J_M(0)$ is computed only once for each fixed model, $d$, and $n$, and can then be used throughout the entire central region
$\|\Sigma^{-1/2}a\|\le L(d/n)^{1/2}$.

\subsection{Local limit theorem and central expansion}\label{ssec:local CLT}
For convenience, in this subsection we suppose that $X$ is isotropic, i.e., $\Eb X=0$ and $\Eb (XX^\top)=I_d$. A standard result on the Legendre transform (which we derive for convenience in Appendix~\ref{ssec:legendre}) implies
\be\label{poly-bound}
\bigl|\varphi^*(a)-\tfrac12\Vert a\Vert^2\bigr|
\les C_3(a)\Vert a\Vert^3,
\ee
where $\varphi^*$ is defined in \eqref{rho a} and $a\in U_d$ satisfies
\be\label{C3 def}
4\Vert a\Vert C_3(a)\le 1,\quad 
C_3(a):=\sup_{\Vert \tau\Vert\le 2\Vert a\Vert}\Vert\nabla^3 \varphi(\tau)\Vert_{\op}.
\ee
We will only consider those $a$ for which all the $\tau$ appearing in the supremum above belong to $V_d$. By \eqref{rho a}, the following local CLT is a corollary to Theorem~\ref{thm:mainI}.

\begin{corollary}\label{cor:CLT} Suppose $X$ is isotropic and the assumptions
of Theorem~\ref{thm:mainI} hold. For any $x\in\br^d$ and $n\in\BN$ such that 
\begin{enumerate}
\item $\{\tau\in\br^d:\|\tau\|\le 2\|x\|/n^{1/2}\}\subset V_d$, and 
\item \eqref{C3 def} is satisfied with $a=x/n^{1/2}\in U_d$, 
\end{enumerate}
one has
\bs\label{rhon}
\rho_n(a)\dd a=&I(a)e^{-\mu(a)}\ga(\dd x),\ a=x/n^{1/2},\\
\mu(a)=&n\big(\varphi^*(a)-\tfrac12\|a\|^2\big)+\tfrac12\log\det H(\tau),
\es
where $|I(a)-1|$ is bounded by \eqref{main res I}, and 
\be\label{mu bnd}
|\mu(a)|\les C_3(a)\frac{\Vert x\Vert^3+d\|x\|}{n^{1/2}}.
\ee
\end{corollary}

\begin{proof} In view of \eqref{poly-bound}, it only remains to show that the determinant in \eqref{rho a} gives the second term in the fraction in \eqref{mu bnd}. Since $H(0)=I_d$, we have
\be
H(\tau)-I_d =\int_0^1 \nabla^3\varphi(s\tau)[\tau,\cdot,\cdot]\dd s,
\ee
and hence, using $\|\tau\|\le 2\|a\|$,
\be
\|H(\tau)-I_d\|\le C_3(a)\|\tau\| \le 2C_3(a)\|a\|\le \tfrac12.
\ee
Thus, all the eigenvalues of $H(\tau)$ are bounded away from zero and
\be
|\log\det H(\tau)| \les d\|H(\tau)-I_d\| \les d C_3(a)\|a\|.
\ee
Setting $a=x/\sqrt n$ completes the proof. 
\end{proof}

The leading central approximation to the density is, of course, the Gaussian approximation furnished by the CLT.  A natural next step is to retain the first non-Gaussian terms in the central expansion of the saddlepoint representation \eqref{rhon}. These terms are determined by the third cumulant tensor and give a cubic correction to the Gaussian
approximation, as shown in the following corollary.

\begin{corollary}\label{cor:central density expansion}
Suppose $X$ is isotropic and \eqref{uniform directional moments} holds for some $r>0$ independent of $d$. Suppose also that Assumption~\ref{ass:mainI}(\ref{global}) holds with $V_d=\{\|\tau\|<r\}$. Let $T_3:=\nabla^3\varphi(0)$, $f_n$ denote the density of
\be
S_n:=\frac1{\sqrt n}\sum_{j=1}^nX_j,
\ee
and let $\ga_d$ denote the standard Gaussian density on $\br^d$. Then, for every fixed $L>0$, uniformly for $\|x\|\le L\sqrt d$, one has
\be\label{central density expansion}
f_n(x)
=
\ga_d(x)
\exp\big\{
\tfrac1{6\sqrt n}T_3[x,x,x]
-\tfrac1{2\sqrt n}\tr T_3[x,\cdot,\cdot]
\big\}
\big(1+O_L(d^2/n)+\est\big),
\ee
provided $d^2/n$ is sufficiently small.  Here
\be\label{third tensor trace}
\tr T_3[x,\cdot,\cdot]
:=
\sum_{j=1}^d T_3[x,e_j,e_j].
\ee
\end{corollary}

The proof of the corollary is straightforward and is given in Appendix~\ref{sec:prf cor expansion}. 

\subsection{More general SPA theorem}\label{sec:resII}

Let $Y_n$, $n\ge 1$, be a sequence of random vectors in $\br^d$, where $d$ may grow with $n$. Denote by
\be\label{GFs II}
\phi_n(\tau)=\big(\Eb e^{n\langle \tau, Y_n\rangle}\big)^{1/n},\quad \varphi_n(\tau)=\log\phi_n(\tau).
\ee
For complex arguments, the powers in \eqref{GFs II} are understood via the analytic branch obtained by continuation from the positive real axis. The random vector $Y_n$ is an analog of the average $n^{-1}\sum_{i=1}^n X_i$ in Section~\ref{sec:setting}. We suppose that $\phi_n(t)$, $n\ge 1$, satisfy the following assumption. 

\begin{assumption}\label{ass:mainII} For all $n$ sufficiently large,
\begin{enumerate}
\item $\sup_{\tau\in V_d}\|\nabla^k \phi_n(\tau)\|_{\op}<\infty$, $0\le k\le4$, for some convex neighborhood of the origin $0\in V_d\subset\br^d$.
\item $H_n(\tau):=\nabla^2\varphi_n(\tau)>0$ for all $\tau\in V_d$ and $\sup_{\tau\in V_d}\Vert H_n^{-1}(\tau)\Vert<\infty$.
\item\label{globalII} One has
\bs\label{ass_1II}
\bigg|\frac{\phi_n(\tau+iH_n^{-1/2}  t)}{\phi_n(\tau)}\bigg|\le & \max\Big(\exp\big\{-\e^{1/2}\Vert t\Vert\big\},(1+\Vert t\Vert)^{-\frac{\kappa}{d^{1/2}}}\Big),\\ 
\tau\in & V_d,\ \Vert t\Vert\ge 2.5\e^{1/2},
\es
for some $\kappa>0$, where $H_n=H_n(\tau)$.
\item\label{phase_v2} There exists $\de>0$ such that 
\be
|\text{Arg}\,\phi_n(\tau+iH_n^{-1/2}  t)|\le\pi-\de,\quad 
|\phi_n(\tau+iH_n^{-1/2}  t)|\ge\de, 
\ee
for all $\tau\in V_d$, $\Vert t\Vert\le 2.5\e^{1/2}$.
\end{enumerate}
\end{assumption}
The constants $\kappa$ and $\de$ are independent of $d$ and $n$, and the size of $V_d$ may depend on $d$. Similarly to \eqref{c3c4}, we define
\bs\label{c3c4 v2}
c_3:=\sup_{\tau\in V_d,\Vert t\Vert < 2.5\e^{1/2}}\Vert \nabla_t^3\im\varphi_n(\tau+i H_n^{-1/2}  t)\Vert_{\op},\\
c_4:=\sup_{\tau\in V_d,\Vert t\Vert < 2.5\e^{1/2}}\Vert \nabla_t^4\re\varphi_n(\tau+i H_n^{-1/2}  t)\Vert_{\op}.
\es
Assumption~\ref{ass:mainII} implies that the constants $c_k$ are finite, but they may grow with $d$ and $n$.

The proof of the following result is identical to that of Theorem~\ref{thm:mainI}, upon replacing $\phi$, $\varphi$, and $H$ with $\phi_n$, $\varphi_n$, and $H_n$, respectively. The quantities $\la_*,A_0,B_0,D_0$ are defined as in \eqref{lots of notation}, where $c_4$ is given in \eqref{c3c4 v2}.  

\begin{theorem}\label{thm:mainII} 
Suppose Assumption~\ref{ass:mainII} is satisfied and $\la_*>0$. For all
$d^2/n$ sufficiently small, one has
\bs\label{main res II}
|I(a)-1|\le &
e^{D_0}\bigg(D_0+\frac{c_3^2}{8\la_*^3}\frac{d(d+4)}n\bigg)
+c\left[\exp(-d)+\bigg(\frac e{\kappa^2}\frac{d^2}n\bigg)^{d/2}\right],\quad
a\in U_d,
\es
for an absolute constant $c>0$.
\end{theorem}

Theorem~\ref{thm:mainII} generalizes Theorem~\ref{thm:mainI} in that it imposes no specific form on $Y_n$; in particular, $Y_n$ need not be the average of $n$ i.i.d. random vectors. A common case is when $Y_n$ is the average of $n$ independent, but not necessarily identically distributed random vectors.

\section{Multivariate generalized Laplace distributions}
\label{sec:VG}

We consider the multivariate generalized Laplace distribution \cite[eq. (4)]{kozubowski2013}. Its pdf is
\bs
\rho(x)=&
\frac{2\exp\{\langle \Sigma^{-1}b,x\rangle\}}
{(2\pi)^{d/2}\Gamma(s)(\det\Sigma)^{1/2}}
\left(
\frac{\langle x,\Sigma^{-1}x\rangle}
     {2+\langle b,\Sigma^{-1}b\rangle}
\right)^{\frac12(s-d/2)}
K_{s-d/2}\left(
\sqrt{(2+\langle b,\Sigma^{-1}b\rangle)
      \langle x,\Sigma^{-1}x\rangle}
\right),
\es
where $b\in\br^d$, $\Sigma>0$, and $K_\nu$ is the modified Bessel function of the second kind. Its mgf and cgf are \cite[eq. (3)]{kozubowski2013}
\begin{align}\label{VG mgf}
\phi(z)=&
\left(1-\langle b,z\rangle-\tfrac12\langle z,\Sigma z\rangle
\right)^{-s},\\
\label{VG cgf}
\varphi(z)
=& -s\log\left(1-\langle b,z\rangle-\tfrac12\langle z,\Sigma z\rangle\right).
\end{align}

We suppose that, uniformly in $d$,
\be\label{VG param bnds}
 \|b\|\les 1,\quad
 0<\sigma_-\le\lambda_{\min}(\Sigma)
 \le\lambda_{\max}(\Sigma)\le\sigma_+<\infty.
\ee
For simplicity, in what follows we assume that $s>0$ is fixed. 

The density has a singularity at the origin once $d>2s$. Indeed, setting $q(x):=\langle x,\Sigma^{-1}x\rangle$, and using $K_\nu(t)=K_{-\nu}(t)$ together with
\be
K_\nu(t)\sim 2^{\nu-1}\Gamma(\nu)t^{-\nu},
\quad t\to0^+,\quad \nu>0,
\ee
we obtain, for $d>2s$,
\be
\rho(x)\sim q(x)^{s-d/2}\asymp \|x\|^{2s-d},
\quad x\to0.
\ee
Thus $\rho(x)\to\infty$ as $x\to0$. In particular, in anticipation of Section~\ref{sec:logconcave}, the generalized Laplace densities considered here are not LC whenever $d>2s$, since a full-dimensional LC probability density is bounded.

Choose $r>0$ sufficiently small and set
\be
 V_d=\{\tau\in\br^d:\|\tau\|<r\}.
\ee
Shrinking $r$ if needed, we have
\be\label{VG A bnd}
 A(\tau):=
 1-\langle b,\tau\rangle-\tfrac12\langle\tau,\Sigma\tau\rangle
 \asymp 1,\quad \tau\in V_d,
\ee
with constants independent of $d$.

Differentiating \eqref{VG cgf} gives
\be\label{VG grad Hess}
\nabla\varphi(\tau)= s\frac{b+\Sigma\tau}{A(\tau)},
\quad  
H(\tau) = s\left\{ \frac{\Sigma}{A(\tau)}
 + \frac{(b+\Sigma\tau)(b+\Sigma\tau)^\top}{A(\tau)^2} \right\}.
\ee
It follows from \eqref{VG param bnds}--\eqref{VG grad Hess} that
\be\label{VG H bnd}
 0<\la_-\le\lambda_{\min}(H(\tau))
 \le\lambda_{\max}(H(\tau))\le \la_+<\infty,
 \quad \tau\in V_d,
\ee
uniformly in $d$.  Thus Assumption~\ref{ass:mainI}(\ref{posdef}) is
satisfied.  Assumption~\ref{ass:mainI}(\ref{exist}) follows immediately
from \eqref{VG mgf} after reducing $r$, if necessary.

We next verify Assumption~\ref{ass:mainI}\eqref{global}.  For $u\in\br^d$,
\be
 A(\tau+iu)
 = A(\tau)+\tfrac12\langle u,\Sigma u\rangle
 -i\langle b+\Sigma\tau,u\rangle.
\ee
Hence
\be\label{VG Fourier exact}
 \left| \frac{\phi(\tau+iu)}{\phi(\tau)} \right|
 =
\bigg[
 \frac{A(\tau)}
 {\big[\left(A(\tau)+\frac12\langle u,\Sigma u\rangle\right)^2
 + \langle b+\Sigma\tau,u\rangle^2
 \big]^{1/2}} \bigg]^s
 \le
 \left( 1+\frac{\langle u,\Sigma u\rangle}{2A(\tau)} \right)^{-s}.
\ee
Take $u=H(\tau)^{-1/2}t$. By \eqref{VG param bnds}, \eqref{VG A bnd}, and \eqref{VG H bnd},
\be
 \langle H(\tau)^{-1/2}t,\Sigma H(\tau)^{-1/2}t\rangle
 \ge c\|t\|^2.
\ee
Therefore
\be\label{VG Fourier bnd}
\Big|\frac{\phi(\tau+iH(\tau)^{-1/2}t)}{\phi(\tau)}\Big|
\le(1+c\|t\|^2)^{-s}.
\ee
For $\|t\|$ bounded away from zero, \eqref{VG Fourier bnd} is much
stronger than the second bound in \eqref{ass_1}.  For small $t$, the
Hessian normalization gives as in \eqref{G t-exp}, uniformly in $\tau$,
\be
\log\Big| \frac{\phi(\tau+iH(\tau)^{-1/2}t)}{\phi(\tau)} \Big|
= -\tfrac12\|t\|^2+O(\|t\|^4).
\ee
Consequently, for some $r_1>0$,
\be\label{VG local Fourier}
\Big| \frac{\phi(\tau+iH(\tau)^{-1/2}t)}{\phi(\tau)} \Big|
\le e^{-0.4\|t\|^2},  \quad \|t\|\le r_1.
\ee
If $2.5\e^{1/2}\le\|t\|\le r_1$, then
\be
 0.4\|t\|^2\ge\e^{1/2}\|t\|.
\ee
Combining this observation with \eqref{VG Fourier bnd} verifies
Assumption~\ref{ass:mainI}(\ref{global}) for all $\e$ sufficiently small.

The quantities $c_3$ and $c_4$ are uniformly bounded as well.  Indeed,
all derivatives of order three and four of \eqref{VG cgf} are finite
sums of terms involving $\Sigma$, $b+\Sigma z$, and negative powers of
$A(z)$.  In the region
\be
 z=\tau+iH(\tau)^{-1/2}t,\quad
 \tau\in V_d,\quad \|t\|\le2.5\e^{1/2},
\ee
the quantities $\|b+\Sigma z\|$, $\|\Sigma\|$, and $|A(z)|^{-1}$ are
uniformly bounded.  Thus, using \eqref{VG H bnd},
\be\label{VG c34}
 c_3+c_4\les1.
\ee
Finally, since $A(\tau+iu)$ remains in a fixed neighborhood of the
positive real axis when $\|u\|\les\e^{1/2}$, the nonvanishing and
phase requirements in Assumption~\ref{ass:mainI}(\ref{phase}) follow
for all $\e$ sufficiently small.

We have therefore proved the following.

\begin{proposition}\label{prop:VG SPA}
Suppose \eqref{VG param bnds} holds and let $X$ have the mgf \eqref{VG mgf}.  Then, setting $V_d=\{\|\tau\|< r\}$ for a sufficiently small fixed $r>0$,  Assumption~\ref{ass:mainI} holds. Moreover, $c_3,c_4$ are bounded independently of $d$ and $n$ and, provided $d^2/n$ is sufficiently small, 
\be\label{VG SPA conclusion}
 |I(a)-1|\les\tfrac{d^2}{n}+\est,
 \quad a\in\nabla\varphi(V_d).
\ee
\end{proposition}

Let us now compare \eqref{VG SPA conclusion} with the high-dimensional saddlepoint theorem of \cite{tang2025}. To map our notation to theirs, let
\be
\mathcal X_n:=\sum_{j=1}^n X_j,\qquad p=d,\qquad s_n=na.
\ee
Then the cgf of $\mathcal X_n$ is $K_{\mathcal X_n}(z)=n\varphi(z)$, and the saddlepoint $\widehat t_n$ in \cite{tang2025} is precisely our
$\tau$, since
\be
\nabla K_{\mathcal X_n}(\tau)=n\nabla\varphi(\tau)=na=s_n.
\ee
Writing, as in \cite{tang2025},
\be
K_{\mathcal X_n}(x+iy)=U(x,y)+iV(x,y),
\ee
and setting
\be
H_r(z):=\re\nabla^2\varphi(z),
\ee
we therefore have
\be
U^{(x,2)}(\tau,y)
=
n\re\nabla^2\varphi(\tau+iy)
=
nH_r(\tau+iy).
\ee
Thus the Hessian appearing in \cite[Assumption~8]{tang2025} is the same as $nH_r(z)$ in our notation.

The bounds \eqref{VG H bnd} and \eqref{VG c34}, together with the same argument on a sufficiently small fixed complex neighborhood, show that the second, third, and fourth derivatives appearing in \cite[Assumptions~8--10]{tang2025} have the required spectral bounds. Moreover, \eqref{VG Fourier bnd}, after raising its right-hand side to the $n$th power, gives the required estimate of the outer Fourier integral in their Assumption~7.

For a matrix $A=(A_{ij})$, we denote by $\|A\|_\infty$ the induced
$\ell_\infty$ operator norm,
\be
\|A\|_\infty
:=
\sup_{\|x\|_\infty=1}\|Ax\|_\infty
=
\max_{1\le i\le d}\sum_{j=1}^d |A_{ij}|.
\ee
After reducing the fixed complex neighborhood if necessary, the argument above \eqref{VG H bnd} gives uniformly $H_r(z)\asymp I_d$. Consequently,
\be\label{TR infinity bound}
\|(nH_r(z))^{-1/2}\|_\infty
\le d^{1/2}\|(nH_r(z))^{-1/2}\|
\les d^{1/2}n^{-1/2}.
\ee
The factor $d^{1/2}$ in \eqref{TR infinity bound} cannot in general be removed under \eqref{VG param bnds}. Indeed, for admissible choices of $b$ and $\Sigma$ one has
\be
\|(nH_r(0))^{-1/2}\|_\infty
\gtrsim (d/n)^{1/2};
\ee
see Appendix~\ref{ssec:TR infinity lower bound}. Thus, in the notation of \cite{tang2025}, one should in general take $c_\infty=1/2$.

Next, consider Assumptions~9--10 in \cite{tang2025}. To avoid confusion with the constants $c_3,c_4$ in \eqref{c3c4}, denote the exponents appearing in these assumptions by $c_3^{\rm TR}$ and $c_4^{\rm TR}$. On the fixed complex neighborhood considered above, direct differentiation of \eqref{VG cgf} gives
\be
\|\nabla^3\varphi(z)\|_{\op}
+
\|\nabla^4\varphi(z)\|_{\op}
\les 1.
\ee
Since $K_{\mathcal X_n}(z)=n\varphi(z)$,
\be
\|\nabla^kK_{\mathcal X_n}(z)\|_{\op}\les n,
\qquad k=3,4.
\ee
Each matrix slice of either derivative tensor has spectral norm bounded by the operator norm of the corresponding tensor. Hence the matrix slices of $U^{(x,3)}$, $U^{(x,4)}$, and $V^{(x,4)}$ appearing in \cite[Assumptions~9--10]{tang2025} all have spectral norm $\lesssim n$. Thus one may take $c_3^{\rm TR}=c_4^{\rm TR}=1$. These are the smallest exponents that can be guaranteed uniformly under \eqref{VG param bnds}: admissible choices of $b$ and $\Sigma$ give nonzero third- and fourth-order derivatives of $\varphi$ of order one, and hence corresponding derivatives of $n\re\varphi$ of order $n$. 

Thus, with $c_\infty=1/2$ and $c_3^{\rm TR}=c_4^{\rm TR}=1$, the error in \cite[Theorem~5.1]{tang2025} becomes
\be
|I(a)-1|
\les
\max\left\{\frac{d^4}{n},\frac{d^3}{n}\right\}
\les \frac{d^4}{n}.
\ee

By contrast, \eqref{VG SPA conclusion} gives the error $O(d^2/n)$ under the same parameter bounds \eqref{VG param bnds}, and hence convergence under the substantially weaker condition $d^2/n\to0$. If additional structure is imposed so that the $\ell_\infty$ Hessian bound in \cite{tang2025} holds with $c_\infty=0$, their error improves to $O(d^3/n)$, which is still one power of $d$ larger than \eqref{VG SPA conclusion}.

\section{A sufficient criterion}\label{sec:sufficient criterion}

Let $\rho(x)$ be the density of $X$. For $\Theta\in S^{d-1}$ and $\tau\in V_d$, define the 1D tilted marginal
\bs\label{two rhos}
\rho_\Theta^\tau(s):=&\tfrac1{\phi(\tau)}\int_{\br^d} e^{\langle \tau,x\rangle}\rho(x)\delta(\langle \Theta,x\rangle-s)\dd x.
\es
Here $\phi(\tau)$ is the normalizing constant. Let $X_\Theta^\tau$ be the random variable with the pdf $\rho_\Theta^\tau$. The quantile $m_{\Theta}^\tau(a)$ of the random variable $\rho_\Theta^\tau(X_\Theta^\tau)$ is any number that satisfies
\be\label{quant def}
\BP\{\rho_\Theta^\tau(X_\Theta^\tau)<m_{\Theta}^\tau(a)\}\le a,\quad \BP\{\rho_\Theta^\tau(X_\Theta^\tau)>m_{\Theta}^\tau(a)\}\le 1-a,\quad 0<a<1.
\ee
\begin{lemma}\label{lem:easy conds}
Suppose \eqref{uniform directional moments} and \eqref{covariance uniform bounds} hold for some $r>0$. Suppose further that
\begin{align}
\label{Holder cond}
&\sup_{\|\tau\|\le r}\sup_{\Theta\in S^{d-1}}
\int_{\br}|\rho_\Theta^\tau(s+h)-\rho_\Theta^\tau(s)|\dd s\les |h|^\al,\quad |h|\le h_0,
\end{align}
for some $0<\al\le 1$ and $h_0>0$, with all constants independent of $d$ and $n$. 
%
Then there exists $0<r'\le r$ such that, with $V_d=\{\|\tau\|<r'\}$, Assumption~\ref{ass:mainI}\eqref{global} holds on $V_d$ for all $\e=d/n$ sufficiently small.
\end{lemma}


\begin{proof}[Proof of Lemma~\ref{lem:easy conds}] 
Choose $r'\in(0,r]$ sufficiently small so that \eqref{spectral bounds} holds. From the assumptions of the lemma and \eqref{cosh bound} it follows that
\be\label{tilted2moments}
\Eb |X_\Theta^\tau|^j\les 1,\quad 1\le j\le 4,\ \|\tau\|<r',\ \Theta\in S^{d-1}.
\ee
Next, consider quantiles of $\rho_\Theta^\tau$. By an elementary calculation, \eqref{Holder cond} yields
\be\label{FT bnd}
|\widetilde{\rho_\Theta^\tau}(\la)|\les (1+|\la|)^{-\al},\quad \|\tau\|<r,
\ee
uniformly in $\tau$. Suppose first that $0<\al<1$. The Besov embedding $B_{1,\infty}^\al(\br)\hookrightarrow L^{q,\infty}(\br)$, where $q=1/(1-\al)>1$ \cite[Part (iii) of Theorem 1.15]{haroske2009}, implies that 
\be\label{Lbnd}
L(m):=|\{s:\rho_\Theta^\tau(s)>m\}|\les m^{-q},\quad m\to\infty.
\ee
Applying this in \eqref{quant def} we get (with $\rho_\Theta^\tau$ denoted by $\rho$ for simplicity)
\be
\BP\{\rho(X)>m\}=\int_{\rho(s)>m}\rho(s)\dd s=mL(m)+\int_m^\infty L(u)\dd u.
\ee
Using \eqref{Lbnd} gives an upper bound 
\be
\BP\{\rho(X)>m\} \les m^{1-q}+\int_m^\infty u^{-q}\dd u\to0,\quad m\to\infty.
\ee
Thus, $m_{\Theta}^\tau(1/2) \les 1$ uniformly in $\Theta\in S^{d-1}$ and $\|\tau\|\le r$.

If $\alpha=1$, then \eqref{Holder cond} implies a uniform $BV$ bound on $\rho^\tau_\Theta$, and hence $\|\rho^\tau_\Theta\|_\infty\lesssim1$; in particular,
$m^\tau_\Theta(1/2)\lesssim1$ directly.

Comparing \eqref{two rhos} with \eqref{ass_1}, we see that the ratio in \eqref{ass_1} is precisely $|\widetilde{\rho_\Theta^\tau}(\la)|$, where $\Theta$ and $\la$ satisfy
$H(\tau)^{-1/2}t=\la\Theta$ and
\be\label{la t}
\Var(X_\Theta^\tau)\la^2=\|t\|^2.
\ee
This relation between $\la$ and $t$ is assumed in what follows. 

To prove Assumption~\ref{ass:mainI}\eqref{global}, begin by considering the interval $2.5\e^{1/2}\le \|t\|\le t_0$ for a sufficiently small $t_0>0$. The function $|\widetilde{\rho_\Theta^\tau}|$ is even, because $\rho_\Theta^\tau$ is real-valued. Using \eqref{tilted2moments} we get 
\be
|\widetilde{\rho_\Theta^\tau}(\la)| \le 1-\tfrac12 \Var(X_\Theta^\tau)\la^2+O(\la^4),\quad \la\to0. 
\ee
 By \eqref{spectral bounds}, $\Var(X_\Theta^\tau)=\langle H(\tau)\Theta,\Theta\rangle\asymp 1$ and hence $\la\asymp\|t\|$. This gives the inequality
\be
|\widetilde{\rho_\Theta^\tau}(\la)|\le \exp\{-\e^{1/2}\|t\|\},\quad
2.5\e^{1/2}\le \|t\| \le t_0,\ \|\tau\|\le r',\ \Theta\in S^{d-1},
\ee
for sufficiently small $t_0>0$ and $r'>0$.

Next, consider the set $\|t\|\ge t_0$. In view of \eqref{tilted2moments} and the uniform bound on $m_{\Theta}^\tau(1/2)$, \cite[Theorem 2]{bobkov2012} implies 
\be\label{bcg bnd}
|\widetilde{\rho_\Theta^\tau}(\la)|\le 
\begin{cases}
1-b_1\la^2,& |\la|\le b_3,\\
1-b_2,& |\la|\ge b_3,
\end{cases}
\ee
for some $b_k>0$, $k=1,2,3$. Using the bound \eqref{FT bnd}, we can find $t_1\ge t_0$, independent of $d$ and $n$, such that
\be
|\widetilde{\rho_\Theta^\tau}(\la)|\les (1+|\la|)^{-\al} \le (1+\|t\|)^{-\kappa/d^{1/2}},\quad \|t\|\ge t_1,\ \|\tau\|\le r',\ \Theta\in S^{d-1},
\ee
for all $d$ sufficiently large. 

Finally, by \eqref{bcg bnd}, \eqref{ass_1} holds for $t_0\le \|t\|\le t_1$ for all $\e$ sufficiently small and $d$ sufficiently large. 
\end{proof}

\section{Structured classes of distributions}
\label{sec:logconcave}

We now verify the dimension-uniform conditions of Lemma~\ref{lem:easy conds} for several structured classes of distributions.

\subsection{LC isotropic case}\label{ssec:LC iso}

Suppose that $X$ is isotropic and the pdf of $X$ is $e^{-V(x)}$, where $V:\br^d\to (-\infty,\infty]$ is extended value convex. Allowing $V$ to take the value $+\infty$ permits, in particular, compactly supported LC densities. Since $X$ is isotropic, \eqref{covariance uniform bounds} trivially holds with absolute constants.\\

\noindent
{\bf Proof of \eqref{uniform directional moments}.} Using that marginals of a LC isotropic random vector are LC isotropic, it suffices to prove \eqref{uniform directional moments} only in $\br$, i.e., for its one-dimensional marginal $e^{-V_1(x)}$, $x\in\br$, on the subspace spanned by $\tau$. 

In $\br$, \cite[Theorem 5.14]{LV2007} implies that $V_1(0)\le 7\log 2$. Consider the sublevel set $D:=\{x:V_1(x)\le V_1(0)+1\}$, which is an interval. The bound on $V_1(0)$ implies that $|D|\le R$, where $|D|$ is the length of the interval and $R>0$ is an absolute constant. Since $V_1$ is a convex function, this means that the rate of linear growth of $V_1$ outside $D$ is bounded below by some absolute positive constant, say, $\mu$. Therefore,
\be\label{V lower bnd}
V_1(x)\ge V_1(0)+1+\mu\text{dist}(x,D),\quad x\not\in D.
\ee
Selecting $\de_2<0.9\mu$ gives
\be\label{1d est}
\int_{\br}e^{\tau x-V_1(x)}\dd x\les 1,\quad |\tau|\le\de_2,
\ee
where the absorbed bound is absolute. This proves \eqref{uniform directional moments} with an absolute bound.

In \eqref{1d est} and similarly in the proof below, we tacitly identify $\tau\in\br^d$ with its projection to the appropriate one- or two-dimensional subspace.
\\

\noindent
{\bf Proof of \eqref{Holder cond}.} To prove \eqref{Holder cond} it suffices to consider the two-dimensional marginal $e^{-V_2(x)}$, $x\in\br^2$, on the subspace spanned by $\Theta$ and $\tau$. If the two are linearly dependent, we can take any subspace containing them. 

In $\br^2$, \cite[Theorem 5.14 and Lemma 5.18]{LV2007} imply that $\|x_0\|\le 3$ and $V_2(0)\le 14\log 2$, where $x_0$ is the location of the minimum of $V_2(x)$. Consider the sublevel set $D:=\{x:V_2(x)\le V_2(x_0)+1\}$. Let $s\in(0,1]$ be its probability. Given that $V_2(x_0)\le V_2(0)$, the bound on $V_2(0)$ implies that $s\ge c_1 |D|$, where $|D|$ is the area of $D$ and $c_1>0$ is an absolute constant. 

By \cite[Lemma 5.13]{LV2007}, $e^{-V_2(x)}$ is $1/e$-rounded, so $D$ contains a ball of radius $s/e$. Pick any $x'\in D$, and let $x_c$ be the center of the ball. Denote $R'=\|x_c-x'\|$. Clearly, $|D|\ge c_2 sR'$ for an absolute constant $c_2>0$. To see this, we consider the triangle formed by $x'$ and the diameter of the ball orthogonal to the line through $x_c$ and $x'$. The triangle is inside $D$, because $D$ is convex. Thus, $s\ge c_1 c_2sR'$, so $R'$ is bounded by an absolute constant, i.e., the diameter of $D$ is bounded by an absolute constant. 

In turn, this implies that there exists an absolute constant $R$ such that $D\subset\{\|x-x_0\|\le R\}$. Since $V_2$ is a convex function, this means that the rate of linear growth of $V_2$ along any ray $x_0+t\Theta$, $t\ge R$, is bounded below by some absolute positive constant, say, $\mu$. Selecting $\de_2<0.9\mu$ we get similarly to the 1D case:
\be\label{exp bnd 2d}
e^{\langle\tau,x\rangle-V_2(x)}\le c_1 e^{-c_2\|x-x_0\|},\quad x\in\br^2,\ \|\tau\|\le\de_2,
\ee
for some absolute constants $c_1,c_2>0$. Furthermore, by the Jensen inequality applied to the convex function $\psi(x)=e^{\langle\tau,x\rangle}$, the normalization constant satisfies
\be
\int_{\br^2} e^{\langle\tau,x\rangle-V_2(x)}\dd x
=\Eb e^{\langle\tau,X\rangle}\ge e^{\langle\tau,\Eb X\rangle}=1,
\ee
where we have used that $\Eb X=0$. Integrating the bound in \eqref{exp bnd 2d} along any line and using the last inequality proves that $\rho_\Theta^\tau(s)$ is bounded by an absolute constant. 

Since the original $d$-dimensional density is LC, its pdf tilted by $e^{\langle\tau,X\rangle}$ is still LC, but it is no longer isotropic. Thus, $\rho_\Theta^\tau$ is also LC (and not necessarily isotropic). Let $s_0$ be a location of its maximum (this point is not necessarily unique). The fact that $\rho_\Theta^\tau$ is LC implies that it is nondecreasing on $(-\infty,s_0]$ and nonincreasing on $[s_0,\infty)$. Therefore,
\bs
\label{Holder cond v3}
\int_{\br}|\rho_\Theta^\tau(s+h)-\rho_\Theta^\tau(s)|\dd s
\le & 2\int_{|s-s_0|\le |h|}\rho_\Theta^\tau(s)\dd s.
\es
Using the uniform bound on $\rho^\tau_\Theta$ obtained above, the right-hand side is $\les |h|$. Hence \eqref{Holder cond} holds with $\al=1$ (taking, for example, $h_0=1$).


Thus, any isotropic LC random variable $X$ satisfies the assumptions of Lemma~\ref{lem:easy conds}. We summarize the result.

\begin{proposition}\label{prop:logconcave SPA}
Let $X=X_d$ be an isotropic random vector with a full-dimensional LC density on $\br^d$. There exists a universal $r_*>0$ such that \eqref{uniform directional moments} and \eqref{Holder cond} hold with $r=r_*$ and $\al=1$, and with $V_d=\{\|\tau\|<r_*\}$ one has
\be\label{LC SPA conclusion}
|I(a)-1| \les\tfrac{d^2}{n}+\est, \quad a\in\nabla\varphi(V_d),
\ee
provided $d^2/n$ is sufficiently small. Finally, $c_3,c_4$ are bounded by a universal constant.
\end{proposition}

\begin{proof}
Since $X$ is isotropic, \eqref{covariance uniform bounds} clearly holds with absolute constants. The preceding argument verifies the remaining assumptions of Lemmas~\ref{lem:aux ass} and \ref{lem:easy conds} with absolute constants (for example, taking $h_0=1$ in \eqref{Holder cond}). After decreasing the resulting universal radius if necessary, denote it by $r_*$. Lemmas~\ref{lem:aux ass} and~\ref{lem:easy conds}, together with Theorem~\ref{thm:mainI}, then prove the proposition.
\end{proof}

\subsection{Affine images and convex mixtures of LC random vectors}
\label{ssec:affine mixtures}

The preceding result is stable under uniformly bounded affine transformations.  This observation also gives a broad class of non-LC distributions to which Lemma~\ref{lem:easy conds} applies.

\begin{lemma}\label{lem:affine logconcave}
Let $X=X_d$ have a full-dimensional LC density on $\br^d$, and suppose that
\be\label{affine LC bnd}
 \|\E X\|\le c_0,\quad
 c_1 I_d\preceq \Var(X)\preceq c_2 I_d,
\ee
with $c_k>0$ independent of $d$. Then, for $\al=1$, there exist $r,h_0>0$ such that \eqref{uniform directional moments} and \eqref{Holder cond} hold. Consequently, the conclusions of Proposition~\ref{prop:logconcave SPA} hold for $X$, with constants depending only on $c_0,c_1,c_2$.
\end{lemma}

\begin{proof}
Set
\be
b=\Eb X,\quad A=\Var(X)^{1/2},\quad Y=A^{-1}(X-b).
\ee
Clearly, $Y$ is isotropic and LC. The argument preceding Proposition~\ref{prop:logconcave SPA} shows that there are absolute constants $r_0,h_0^{(Y)}>0$ such that 
\be
 \sup_{\|\sigma\|\le r_0}\Eb e^{\langle\sigma,Y\rangle}\les 1,
\label{affine LC mgf}
\ee
and \eqref{Holder cond} holds for $Y$, uniformly in $d$.

Fix $\Theta\in S^{d-1}$ and set
\be
 r_\Theta:=\|A^\top\Theta\|,\quad
 \widetilde\Theta:=\frac{A^\top\Theta}{r_\Theta},\quad
 \sigma:=A^\top\tau.
\label{affine LC params}
\ee
By \eqref{affine LC bnd},
\be
 c_1^{1/2}\le r_\Theta\le c_2^{1/2},\quad
 \|\sigma\|\le c_2^{1/2}\|\tau\|.
\label{affine LC params bnd}
\ee
Moreover,
\be
 \langle\tau,X\rangle
 =\langle\sigma,Y\rangle+\langle\tau,b\rangle.
\label{affine LC inner}
\ee
Therefore, after taking $r\le r_0/c_2^{1/2}$ and decreasing it if
necessary,
\bs
 \Eb e^{\langle\tau,X\rangle} 
 \les&
 e^{c_0\|\tau\|} \sup_{\|\sigma\|\le r_0}\Eb e^{\langle\sigma,Y\rangle} \les 1,\quad \|\tau\|<r.
\es
Thus \eqref{uniform directional moments} holds uniformly. 

The assumption \eqref{covariance uniform bounds} is the same as the second half of \eqref{affine LC bnd}. It remains to verify \eqref{Holder cond}.  

Write $\rho_{\Theta,X}^\tau$ and $\rho_{\Theta,Y}^\sigma$ for the tilted one-dimensional marginals defined in \eqref{two rhos}.  Exponential tilting of $X$ by $\tau$ corresponds to exponential tilting of $Y$ by $\sigma=A^\top\tau$.  Hence
\be
 \rho_{\Theta,X}^\tau(s)
 =
 \tfrac1{r_\Theta}
 \rho_{\widetilde\Theta,Y}^{\sigma}
 \big( \tfrac{s-\langle\Theta,b\rangle}{r_\Theta} \big).
\label{affine tilted marginal}
\ee
Consequently,
\bs
 \int_{\br}
 |\rho_{\Theta,X}^\tau(s+h)-\rho_{\Theta,X}^\tau(s)|\dd s
 =&
 \int_{\br}
 \big|
 \rho_{\widetilde\Theta,Y}^{\sigma}
 \big(u+\tfrac{h}{r_\Theta}\big)
 -  \rho_{\widetilde\Theta,Y}^{\sigma}(u)
 \big|\dd u
 \les
 \big|\tfrac{h}{r_\Theta}\big|  \les |h|,
\es
uniformly for $\|\tau\|\le r$ and $|h|$ sufficiently small. This proves \eqref{Holder cond}, where we can take $h_0=c_1^{1/2}h_0^{(Y)}$. Here we use that $r_\Theta\ge c_1^{1/2}$; see \eqref{affine LC params bnd}. 

Lemmas~\ref{lem:aux ass} and~\ref{lem:easy conds}, together with Theorem~\ref{thm:mainI}, then prove the lemma.
\end{proof}

We can now pass from individual LC densities to their convex combinations.  The components need not have the same mean or covariance matrix.

\begin{proposition}
\label{prop:logconcave mixtures}
For $1\le j\le J$, let $X_j=X_{j,d}$ have LC density $\rho_j$ and satisfy
\be\label{mixture component bnd}
 \|\E X_j\|\le c_0,\quad
 c_1I_d\preceq\Var(X_j)\preceq c_2I_d,
\ee
for some constants $c_k>0$ independent of $j$ and $d$. Let $X$ have density
\be
 \rho(x)=\sum_{j=1}^J p_j\rho_j(x),\quad
 p_j\ge0,\quad
 \sum_{j=1}^J p_j=1.
\label{LC mixture}
\ee
Then \eqref{covariance uniform bounds} holds. Moreover, for $\al=1$, there exist $r,h_0>0$ such that \eqref{uniform directional moments} and \eqref{Holder cond} hold. Consequently, the conclusions of Proposition~\ref{prop:logconcave SPA} hold for $X$, with constants depending only on $c_0,c_1,c_2$, and, in particular, independent of $d,J$, and the weights $p_j$. 
\end{proposition}

\begin{proof}
Clearly, $X$ need not be LC, so we verify that it satisfies the assumptions of Lemma~\ref{lem:easy conds}. Denote $b_j=\Eb X_j$. Let
\be
 b:=\Eb X=\sum_{j=1}^J p_jb_j,\quad
 \Sigma:=\Var(X)=\sum_{j=1}^J p_j\big[\Var(X_j)+(b_j-b)(b_j-b)^\top\big].
\label{mixture normalization}
\ee
Thus $\|b\|\le c_0$ and
\be
c_1 I_d\preceq\Sigma\preceq (c_2+c_0^2) I_d.
\label{mixture covariance bnd}
\ee
Indeed,
\be
\sum_{j=1}^J p_j(b_j-b)(b_j-b)^\top
=
\sum_{j=1}^J p_j b_jb_j^\top-bb^\top
\preceq
\sum_{j=1}^J p_j b_jb_j^\top
\preceq c_0^2 I_d,
\ee
which implies the second inequality in \eqref{mixture covariance bnd}. Therefore, $X$ satisfies \eqref{covariance uniform bounds}. 

By Lemma~\ref{lem:affine logconcave}, $X_j$ satisfy \eqref{uniform directional moments} and \eqref{Holder cond} uniformly in $j,d$. By linearity of expectation, $X$ satisfies \eqref{uniform directional moments}. It remains to check \eqref{Holder cond}. We have
\bs
\rho_\Theta^\tau(s)=&\tfrac1{\phi(\tau)}\sum_{j=1}^J p_j\int_{\br^d} e^{\langle\tau,x\rangle}\rho_j(x)\de(\langle\Theta,x\rangle-s)\dd x=\tfrac1{\phi(\tau)}\sum_{j=1}^J p_j\phi_j(\tau)(\rho_j)_\Theta^\tau(s),\\
\phi(\tau)=&\sum_{j=1}^J p_j \phi_j(\tau),\quad \phi_j(\tau):=\E e^{\langle\tau,X_j\rangle}.
\es
Therefore, the triangle inequality and individual bounds imply
\bs
\int_{\br}|\rho_\Theta^\tau(s+h)-\rho_\Theta^\tau(s)|\dd s\les
|h|\tfrac1{\phi(\tau)}\sum_{j=1}^J p_j\phi_j(\tau)=|h|.
\es

Lemmas~\ref{lem:aux ass} and~\ref{lem:easy conds}, together with Theorem~\ref{thm:mainI}, then prove the proposition.
\end{proof}

The same proof also applies to general convex mixtures provided the bounds in \eqref{mixture component bnd} hold uniformly over the mixing parameter.

\subsection{Gaussian-smoothed distributions}
\label{ssec:Gaussian smoothing}

We next consider a broad class of distributions obtained by adding
nondegenerate Gaussian noise to an otherwise arbitrary random vector.  Let
\be\label{Gaussian smoothing model}
 X=Y+Z,\quad
 Z\sim \CN(0,\Sigma),\quad
 Y\perp Z.
\ee
No density or full-dimensionality assumption is imposed on $Y$. The Gaussian component ensures that $X$ has a smooth full-dimensional density.

\begin{proposition}
\label{prop:Gaussian smoothing SPA}
Suppose that, for some $r>0$ independent of $d$,
\be\label{Gaussian smoothing assumptions}
 \sup_{\|\tau\|\le r}\Eb e^{\langle\tau,Y\rangle}\les 1,
 \quad
 \Sigma\asymp I_d.
\ee
Then \eqref{covariance uniform bounds} holds. Moreover, there exist $r',h_0>0$ such that \eqref{uniform directional moments} and \eqref{Holder cond} hold with $\al=1$. Consequently, the conclusions of Proposition~\ref{prop:logconcave SPA} hold for $X$, with constants depending only on $r$ and the absorbed constants in \eqref{Gaussian smoothing assumptions}. 
\end{proposition}

\begin{proof}
Denote the mgf of $Y$ by $\phi_Y$.  Independence in \eqref{Gaussian smoothing model} gives
\be
 \phi(\tau)
 =
 \phi_Y(\tau)
 \exp\left\{\tfrac12\langle\Sigma\tau,\tau\rangle\right\}.
\ee
Therefore, \eqref{Gaussian smoothing assumptions} implies
\eqref{uniform directional moments} for $X$.

Moreover, \eqref{Gaussian smoothing assumptions} and \eqref{cosh bound} give $\Eb\langle\Theta,Y\rangle^2\les1$ uniformly in $\Theta\in S^{d-1}$. Thus $\Var(Y)\lesssim I_d$, and hence
\be
\Var(X)=\Var(Y)+\Sigma\asymp I_d.
\ee
This proves \eqref{covariance uniform bounds}.

It remains to verify \eqref{Holder cond}.  For $\|\tau\|\le r$, let
$Y^\tau$ and $Z^\tau$ denote the exponentially tilted versions of $Y$
and $Z$, respectively.  Since $Y$ and $Z$ are independent, exponential
tilting preserves their independence and
\be\label{Gaussian tilted decomposition}
 X^\tau\stackrel{d}=Y^\tau+Z^\tau,
 \quad
 Z^\tau\sim \CN(\Sigma\tau,\Sigma).
\ee
Fix $\Theta\in S^{d-1}$ and denote by $\mu_\Theta^\tau$ the law of
$\langle\Theta,Y^\tau\rangle$.  Also, let $g_\Theta^\tau$ be the pdf of
$\langle\Theta,Z^\tau\rangle$.  Then
\be\label{Gaussian tilted convolution}
\rho^\tau_\Theta(s)
=
\int_{\br} g^\tau_\Theta(s-u)\mu^\tau_\Theta(\dd u),
\quad s\in\br.
\ee
The Gaussian density $g_\Theta^\tau$ has variance
\be
 \sigma_\Theta^2
 =  \langle\Sigma\Theta,\Theta\rangle\asymp1.
\ee
Consequently, elementary transformations give
\bs
 \int_{\br}
 |\rho_\Theta^\tau(s+h)-\rho_\Theta^\tau(s)|\dd s
 &\le
 \int_{\br}
 |g_\Theta^\tau(s+h)-g_\Theta^\tau(s)|\dd s\\
 &\le
 |h|\int_{\br}|(g_\Theta^\tau)'(s)|\dd s
 \les |h|
\es
uniformly in $\|\tau\|\le r$ and $\Theta\in S^{d-1}$.  
Indeed,
\be
 \int_{\br}|(g_\Theta^\tau)'(s)|\dd s
 =
 \sqrt{\frac{2}{\pi}}\,\frac1{\sigma_\Theta}
 \les1.
\ee
This proves \eqref{Holder cond}.

Lemmas~\ref{lem:aux ass} and~\ref{lem:easy conds}, together with Theorem~\ref{thm:mainI}, then prove the proposition.
\end{proof}

A useful special case of Proposition~\ref{prop:Gaussian smoothing SPA} is
the latent-factor model
\be\label{Gaussian smoothing factor}
 X=A_dU+Z,
\ee
where $A_d:\br^{k_d}\to\br^d$ may have rank strictly smaller than $d$.
Suppose $\|A_d\|\les1$ and, for some $r_0>0$,
\be
 \sup_{\|u\|\le r_0}\Eb e^{\langle u,U\rangle}\les1
\ee
uniformly in $d$.  Then, after taking $r>0$ sufficiently small,
\be
 \sup_{\|\tau\|\le r}
 \Eb e^{\langle\tau,A_dU\rangle}
 =
 \sup_{\|\tau\|\le r}
 \Eb e^{\langle A_d^\top\tau,U\rangle}
 \les1,
\ee
so Proposition~\ref{prop:Gaussian smoothing SPA} applies.  Thus the
latent variable may have dimension $k_d<d$, and may even have a singular
or discrete distribution; the uniformly nondegenerate Gaussian noise is
sufficient to provide the regularity required by
Lemma~\ref{lem:easy conds}.

\begin{remark}
\label{rem:general smoothing}
The Gaussian assumption on $Z$ in Proposition~\ref{prop:Gaussian smoothing SPA} is not essential. More generally, $Z$ may be any random vector satisfying \eqref{uniform directional moments}, \eqref{covariance uniform bounds}, and \eqref{Holder cond} uniformly in $d$. Indeed, the first two properties together with the exponential-moment bound on $Y$ in \eqref{Gaussian smoothing assumptions} imply \eqref{uniform directional moments} and \eqref{covariance uniform bounds} for $X=Y+Z$ exactly as in the proof above, while the convolution representation \eqref{Gaussian tilted convolution} and \eqref{Holder cond} for $Z$ imply \eqref{Holder cond} for $X$. In particular, $Z$ may be any of the random vectors covered by Proposition~\ref{prop:logconcave mixtures}.
%
\end{remark}

\section{Application: non-Gaussian noise in averaged images}
\label{sec:image application}

We conclude with a numerical application of our results to likelihood-based inference from averaged images.  The setup is loosely motivated by astronomical image stacking, where averaging many registered images can lead to an approximately Gaussian distribution by the central limit theorem; see, for example, \cite{lewis2009}. We use this astronomical application as a natural high-dimensional test problem for the approximations developed above.  In particular, the two-component Gaussian-mixture noise model below is chosen because it is non-Gaussian and nonseparable while still allowing the exact likelihood of the averaged image to be computed as a finite Gaussian mixture.

Suppose that the same image is observed repeatedly and that the images
have been registered on a common pixel grid.

\subsection{Noise model and exact likelihood}
\label{ssec:image mixture model}

Let $\Lambda_d=\{1,\ldots,m\}^2$, so that $d=m^2$, and let $m_\theta\in\br^d$ denote the underlying image on the patch. The indices are understood modulo $m$ in both coordinate directions. We observe
\be\label{image mixture obs}
Y_j=m_\theta+X_j,
\quad j=1,\ldots,n,
\ee
where the $X_j$ are i.i.d. realizations of a noise vector $X\in\br^d$.

To simulate spatially correlated noise we use the periodic five-point stencil
\bs\label{image mixture G}
(G_\al z)(x)
=&
c_\al\Big[z(x)+\al\sum_{e\in\{(1,0),(-1,0),(0,1),(0,-1)\}}z(x-e)\Big],
\\
c_\al=&(1+4\al^2)^{-1/2},\quad x\in\Lambda_d,
\es
and set $C_\al:=G_\al G_\al^\top$. The Fourier symbol of $G_\al$ is
\be
c_\al\{1+2\al(\cos\om_1+\cos\om_2)\}.
\ee
Thus, for $|\al|<1/4$, $C_\al\asymp I_d$ uniformly in $d$.

Let $Z_j\sim\CN(0,I_d)$, and let $S_j$ be independent random variables taking
the two values
\be\label{image mixture S}
s_-=-\tfrac1{\sqrt3},
\quad
s_+=\sqrt3,
\quad
\BP\{S_j=s_-\}=\tfrac34,
\quad
\BP\{S_j=s_+\}=\tfrac14.
\ee
Thus $\Eb S_j=0$ and $\Eb S_j^2=1$. Fix a smooth spatial pattern
$\mu\in\br^d$ and define
\be\label{image mixture noise}
X_j
=
\sigma\left(
S_j\mu+\sqrt{c_{S_j}}\,G_\al Z_j
\right),
\quad
c_{s_-}=c_-,
\quad
c_{s_+}=c_+.
\ee
Here $\sigma>0$ is fixed and controls the overall noise magnitude, $c_-,c_+>0$ are fixed, and the $S_j$ and $Z_j$ are mutually independent i.i.d. sequences.

The two values of $S_j$ may be viewed as representing two different observation conditions. The first term in \eqref{image mixture noise} represents a random coherent perturbation along the fixed spatial pattern $\mu$, while the second represents spatially correlated Gaussian noise whose strength is modulated by the same latent state $S_j$.

Conditional on $S_j=s_\pm$,
\be\label{image mixture conditional}
X_j\mid S_j=s_\pm
\sim
\CN\left(
\sigma s_\pm\mu,\,
\sigma^2c_\pm C_\al
\right).
\ee
We choose $c_-,c_+$ so that
\be\label{image mixture cbar}
\bar c:=\tfrac34c_-+\tfrac14c_+=1.
\ee
Consequently,
\be\label{image mixture Sigma}
\Eb X_j=0,
\quad
\Sigma:=\Cov(X_j)
=
\sigma^2(C_\al+\mu\mu^\top).
\ee

The cgf of $X_j$ is
\bs
\label{image mixture cgf}
\varphi(\tau)
=
\log\Big(&
\tfrac34
\exp\big\{
\sigma s_-\langle\mu,\tau\rangle
+\tfrac{\sigma^2c_-}{2}\langle\tau,C_\al\tau\rangle\big\}
\\
&+\tfrac14\exp\big\{\sigma s_+\langle\mu,\tau\rangle
+\tfrac{\sigma^2c_+}{2}\langle\tau,C_\al\tau\rangle\big\}
\Big).
\es
The density of $X_j$ is a convex combination of two Gaussian, hence LC, densities, with means $\sigma s_\pm\mu$ and covariance matrices $\sigma^2c_\pm C_\al$. If $\|\mu\|$ is bounded uniformly in $d$ and $|\al|<1/4$, their means and covariance matrices satisfy the uniform bounds in \eqref{mixture component bnd}. Hence Proposition~\ref{prop:logconcave mixtures} implies that the conclusions of Proposition~\ref{prop:logconcave SPA} hold for the present model.

For the numerical experiment we can also evaluate the true likelihood exactly.
Let
\be
K:=\#\{j:S_j=s_+\},
\quad
K\sim{\rm Bin}(n,1/4).
\ee
Conditional on $K=k$, the average $\bar X_n=n^{-1}\sum_{j=1}^nX_j$ is
Gaussian:
\be\label{image mixture conditional average}
\bar X_n\mid K=k
\sim
\CN(m_k,v_kC_\al),
\ee
where
\be\label{image mixture mk vk}
m_k
=
\frac{\sigma}{n}\big[(n-k)s_-+ks_+\big]\mu,
\quad
v_k
=
\frac{\sigma^2}{n^2}\big[(n-k)c_-+kc_+\big].
\ee

By \eqref{image mixture obs}, $\bar Y=m_\theta+\bar X_n$, where the bars denote the average values, so the likelihood of $\theta$ based on the averaged observation $\bar y$ is the pdf of $\bar X_n$ evaluated at $\bar y-m_\theta$. Writing
\be\label{image mixture pik}
 \pi_k=\binom{n}{k}\left(\frac14\right)^k
 \left(\frac34\right)^{n-k},
\ee
and setting $a_\theta:=\bar y-m_\theta$, we obtain the exact likelihood, with the full density normalization retained,
\bs
\label{image mixture exact likelihood}
L_n^{\rm true}(\theta;\bar y)
=&\rho_n(\bar y-m_\theta)=\rho_n(a_\theta)\\
=&\sum_{k=0}^n
\frac{\pi_k}{(2\pi v_k)^{d/2}(\det C_\al)^{1/2}}
\exp\big\{-\tfrac1{2v_k}\langle a_\theta-m_k,\,C_\al^{-1}(a_\theta-m_k)\rangle\big\}.
\es

\subsection{Conventional SPA and its central expansions}
\label{ssec:image mixture SPA}

For the two-point variable in \eqref{image mixture S},
\be\label{image mixture kappa beta}
\kappa_3:=\Eb S_j^3=2/\sqrt3,
\quad
\beta:=\Eb(S_jc_{S_j}).
\ee
For the numerical values used below, $c_-=0.7$ and $c_+=1.9$, so
$\beta=3\sqrt3/10$.

In addition to the conventional SPA, we consider the first two central
approximations discussed in Subsection~\ref{ssec:local CLT}: the leading
Gaussian approximation and its first non-Gaussian, cubic correction.

Let $T_3:=\nabla^3\varphi(0)$ and put
$u_\theta:=\Sigma^{-1}a_\theta$. The logarithm of the conventional SPA
obtained from \eqref{rho a} by replacing $I(a_\theta)$ by $1$ is
\be\label{image mixture conventional SPA}
\ell_n^{\rm S}(\theta;\bar y)
=
\frac d2\log n-\frac d2\log(2\pi)
-\frac12\log\det H(\tau_\theta)
-n\varphi^*(a_\theta),
\quad
\nabla\varphi(\tau_\theta)=a_\theta.
\ee
Expanding \eqref{image mixture conventional SPA} about the central saddlepoint $a_\theta=0$, $\tau_\theta=0$, the constant term is
\be\label{image mixture constant}
C_{n,\Sigma}
:=
\frac d2\log n-\frac d2\log(2\pi)
-\frac12\log\det\Sigma.
\ee
Using the expansions derived in Appendix~\ref{app:image cubic derivation}, retaining the zeroth-order term in $\log\det H$ and the quadratic term in $\varphi^*$ (see \eqref{app image determinant expansion} and \eqref{app image general expansions}, respectively) gives the Gaussian-order approximation
\be\label{image mixture Gaussian}
\ell_n^{\rm G}(\theta;\bar y)
=
C_{n,\Sigma}
-\frac n2\langle a_\theta,\Sigma^{-1}a_\theta\rangle.
\ee
Retaining also the first non-Gaussian terms in the Legendre transform and in the determinant of the Hessian gives the cubic SPA
\be\label{image mixture cubic}
\ell_n^{\rm C}(\theta;\bar y)
=
\ell_n^{\rm G}(\theta;\bar y)
+\frac n6T_3[u_\theta,u_\theta,u_\theta]
-\frac12\langle q,u_\theta\rangle,
\ee
where
\bs\label{q and T3}
T_3[u,u,u]
:=&
\sigma^3
\left[
\kappa_3\langle\mu,u\rangle^3
+
3\beta \langle\mu,u\rangle
\langle u,C_\al u\rangle
\right],\quad u\in\br^d,\\
q
:=&
\sigma^3
\left[
\kappa_3
\langle\mu,\Sigma^{-1}\mu\rangle\mu
+
\beta\left\{\tr(C_\al \Sigma^{-1})\mu+2C_\al \Sigma^{-1}\mu
\right\}\right].
\es
The first correction in \eqref{image mixture cubic} comes from the cubic term in the Legendre transform, and the second from the expansion of the determinant of the Hessian; see Appendix~\ref{app:image cubic derivation} for more details. The normalization constant $C_{n,\Sigma}$ in \eqref{image mixture constant} is retained because below we compare absolute density values. It could be dropped for maximization over $\theta$.

Fix $R>0$ and define the central parameter region
\be\label{image mixture central region}
\Theta_{R,n}(\bar y)
:=
\big\{
\theta: \|\Sigma^{-1/2}a_\theta\| \le R(d/n)^{1/2}\big\}.
\ee
The same central expansion used in the proof of
Corollary~\ref{cor:central density expansion}, namely
\eqref{central saddle expansion}--\eqref{central logdet expansion}, gives
\be\label{image mixture cubic accuracy}
\ell_n^{\rm true}(\theta;\bar y):=\log L_n^{\rm true}(\theta;\bar y)
=
\ell_n^{\rm C}(\theta;\bar y)
+
O_R(d^2/n)
\ee
uniformly for $\theta\in\Theta_{R,n}(\bar y)$. Thus the cubic SPA retains the leading contribution of the third cumulant.

\subsection{Numerical experiment}
\label{ssec:image mixture numerical}

We take
\bs\label{image mixture dimensions}
&m=6, \quad d=36, \quad n=1000,\\
&\al=0.1, \quad \sigma=0.25, \quad c_-=0.7, \quad c_+=1.9.
\es
Notice that \eqref{image mixture cbar} is satisfied, and $d^2/n=1.296$.

We carry out the numerical experiment in a regime where $d^2/n=1.296$ is of order one rather than asymptotically small. There are two reasons for this choice. First, it provides a fairly stringent finite-dimensional stress test of the asymptotic theory. Second, when
$d^2/n$ is made small, the errors of the SPA and of its central expansions become so small that their differences are difficult to illustrate numerically. The present regime makes these differences visible.

\begin{figure}[t]
\centering
\includegraphics[width=.42\textwidth]{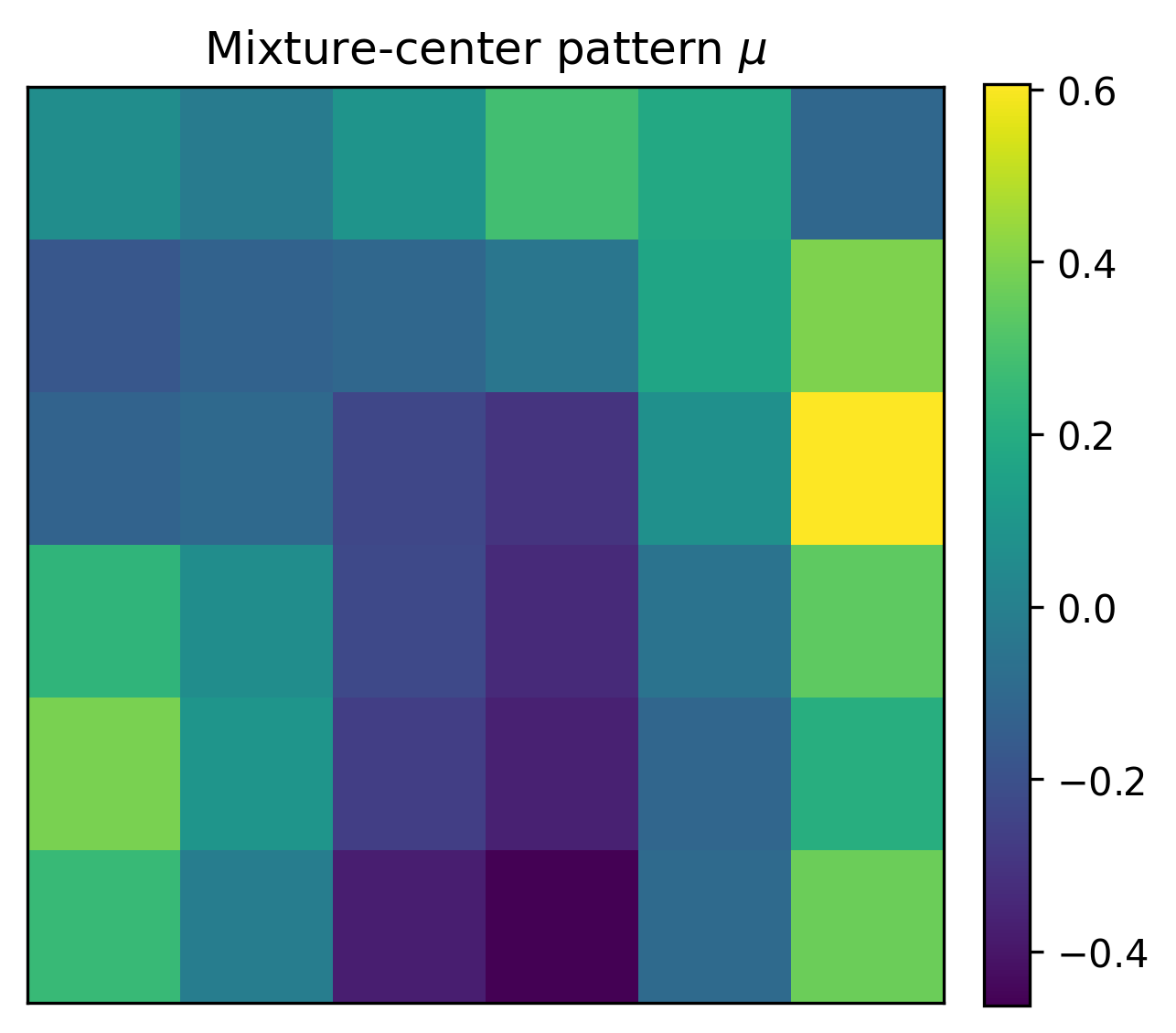}
\caption{The smooth mixture-center pattern $\mu$ used in the numerical
experiment. Its Euclidean norm is $1.5$.}
\label{fig:image mixture mu}
\end{figure}

For later use in the reconstruction experiment, and also to specify the fixed spatial pattern $\mu$, we use the first $4\times4$ tensor-product orthonormal discrete cosine transform (DCT) modes on the $6\times6$ patch:
\be\label{image mixture DCT}
m_\theta=B\theta,
\quad
B:\br^{16}\longrightarrow\br^{36}.
\ee
The vector $\mu$ is a fixed smooth pattern in the same DCT subspace. The 16 DCT coefficients of $\mu$ were generated once with weights proportional to $(1+k+\ell)^{-1/2}$, and the resulting vector was normalized so that $\|\mu\|=1.5$. The resulting pattern is shown in Figure~\ref{fig:image mixture mu}.

The exact likelihood \eqref{image mixture exact likelihood} contains only $n+1=1001$ Gaussian components, and the entire sum is evaluated deterministically. The conventional SPA is evaluated directly from \eqref{rho a} with $I(a)$ replaced by $1$, while the Gaussian and cubic approximations are evaluated from \eqref{image mixture Gaussian} and \eqref{image mixture cubic}.

We first compare the density approximations themselves, independently of any particular ground-truth image or data realization. Write
\be
\rho_n^{\rm S}(a)
:=
\frac{\cdn}{\sqrt{\det H(\tau)}}e^{-n\varphi^*(a)},
\quad
\nabla\varphi(\tau)=a,
\ee
for the conventional SPA obtained from \eqref{rho a} by replacing $I(a)$ by $1$. We similarly write $\rho_n^{\rm G}(a)$ and $\rho_n^{\rm C}(a)$ for the Gaussian and cubic  approximations obtained from \eqref{image mixture Gaussian} and \eqref{image mixture cubic}, respectively, with $a_\theta$ replaced by $a$.

We begin with the locations of the modes. Let $a_*^{\rm true}$, $a_*^{\rm S}$, and $a_*^{\rm G}$ denote the global maximizers of the exact density, the conventional SPA, and the Gaussian approximation, respectively, and let $a_*^{\rm C}$ denote the central local maximizer of the cubic approximation. The latter qualification is necessary because $\rho_n^{\rm C}$ is only a central approximation.

Set
\be\label{e dir def}
e := \frac{\Sigma^{-1/2}\mu}{\|\Sigma^{-1/2}\mu\|}.
\ee
By \eqref{image mixture Sigma} and \eqref{image mixture cgf}, in the whitened coordinates $\rho_n$ is invariant under rotations that fix $e$ (see also \eqref{app image cgf reduced}). The same symmetry is inherited by $\rho_n^M$, where $M\in\{{\rm S},{\rm C},{\rm G}\}$. Writing
\be
\Sigma^{-1/2}a=ze+y,\quad y\perp e,\quad r:=\|y\|,
\ee
this symmetry implies that the exact density and the three approximations
depend only on $(z,r)$. Numerically, the maximizing transverse radius is
$r=0$ in all four cases, so all four relevant modes lie on the axis
spanned by $e$. Define
\be
\zeta_M
:=
\sqrt n\,
\langle e,\Sigma^{-1/2}a_*^M\rangle,
\quad
\delta_M
:= \sqrt n\,
\|\Sigma^{-1/2}(a_*^M-a_*^{\rm true})\|,
\ee
where $M\in\{{\rm true},{\rm S},{\rm C},{\rm G}\}$. Thus $\zeta_M$
measures the signed mode location in units of the one-dimensional
standard deviation of the averaged noise, while $\delta_M$ measures the
error in the mode location. The results are

\begin{table}[ht]
\centering
\begin{tabular}{c|cc}
& $\zeta_M$ & $\delta_M$\\
\hline
Exact density       & $-0.24737$ & $0$\\
Conventional SPA    & $-0.24702$ & $3.48\times10^{-4}$\\
Cubic SPA           & $-0.24630$ & $1.06\times10^{-3}$\\
Gaussian            & $0$        & $2.47\times10^{-1}$
\end{tabular}
\caption{Locations of the exact and approximate modes. For the cubic
approximation, the central local mode is used.}
\label{tab:image mixture modes}
\end{table}

Thus, although $\Eb\bar X_n=0$, the exact density has a nonzero mode as a
result of the skewness of the noise distribution. Both the conventional
SPA and the cubic approximation reproduce this finite-$n$ displacement
very accurately. The Gaussian approximation, whose mode is necessarily
at zero, misses it completely.

We next assess the sharpness of the bound in Theorem~\ref{thm:mainI}. For simplicity, we do this only at the central point $a=0$, for which the corresponding saddlepoint is $\tau=0$. At this fixed saddlepoint, the apparently high-dimensional maximizations in
\eqref{c3c4} reduce, by the rotational symmetry of the whitened model, to four-dimensional compact optimization problems; see Appendix~\ref{app:image c3c4}. Numerically, this gives
\be
c_{3,0}=1.306,\quad c_{4,0}=2.336,
\ee
where the subscript $0$ indicates that the supremum over $\tau\in V_d$ in \eqref{c3c4} is replaced by $\tau=0$. Using these values in the proof of Theorem~\ref{thm:mainI} gives
\be
|I(0)-1|\le 1.064+\est\ .
\ee
Direct numerical evaluation gives instead
\be
|I(0)-1|=1.595\times10^{-2}.
\ee
Thus, disregarding the exponentially small tail terms, the theoretical bound exceeds the actual error at the central point by a factor of approximately $67$. We regard this as a reasonably informative non-asymptotic bound, particularly since it is being evaluated in a regime where $d^2/n=1.296$ is not small.

We next illustrate the Monte Carlo refinement of Subsection~\ref{ssec:MC}. Using $M=10^7$ samples in \eqref{MC estimator}, we obtain
\be
\widehat J_M(0)-1=1.5935\times10^{-2}.
\ee
As an a posteriori validation, comparing with the exact value gives
\be
|I(0)-\widehat J_M(0)|=1.60\times10^{-5}.
\ee
Thus the Monte Carlo calculation gives an accurate estimate of the SPA
error at the central point.

We next use the exact density to compare the conventional SPA and its two central expansions throughout the natural high-dimensional central region. Introduce the standardized shells
\be
\label{image mixture MC shells}
\mathcal S_r
:=
\left\{
a:\|\Sigma^{-1/2}a\|
=
r\big(d/n\big)^{1/2}
\right\}.
\ee
Thus $r$ measures radial distance in units of the natural
high-dimensional central scale $(d/n)^{1/2}$. Under the Gaussian
approximation, the typical radius corresponds to $r\approx1$.

For $M\in\{{\rm S},{\rm C},{\rm G}\}$ define the signed log-density
error
\be\label{image mixture log errors}
\Delta_M(a)
:=
\log\frac{\rho_n^M(a)}{\rho_n(a)},
\ee
and set $\Delta_{\rm true}(a)\equiv0$. In particular,
\be
\Delta_{\rm S}(a)=-\log I(a).
\ee

By the rotational symmetry described above, each $\Delta_M(a)$ depends on $a$ only through $\|\Sigma^{-1/2}a\|$ and $\langle e,\Sigma^{-1/2}a\rangle$. Hence, on each
$\mathcal S_r$, the complete directional dependence is described by the single angular coordinate 
\be\label{image mixture angular coordinate}
\xi:=
\frac{\langle \Sigma^{-1/2}a,e\rangle} {\|\Sigma^{-1/2}a\|} \in[-1,1].
\ee
Here $\xi=1$ corresponds to the direction of $e$, $\xi=-1$ to the opposite direction, and $\xi=0$ to directions orthogonal to $e$. 

\begin{figure}[h]
\centering
\includegraphics[width=.75\textwidth]{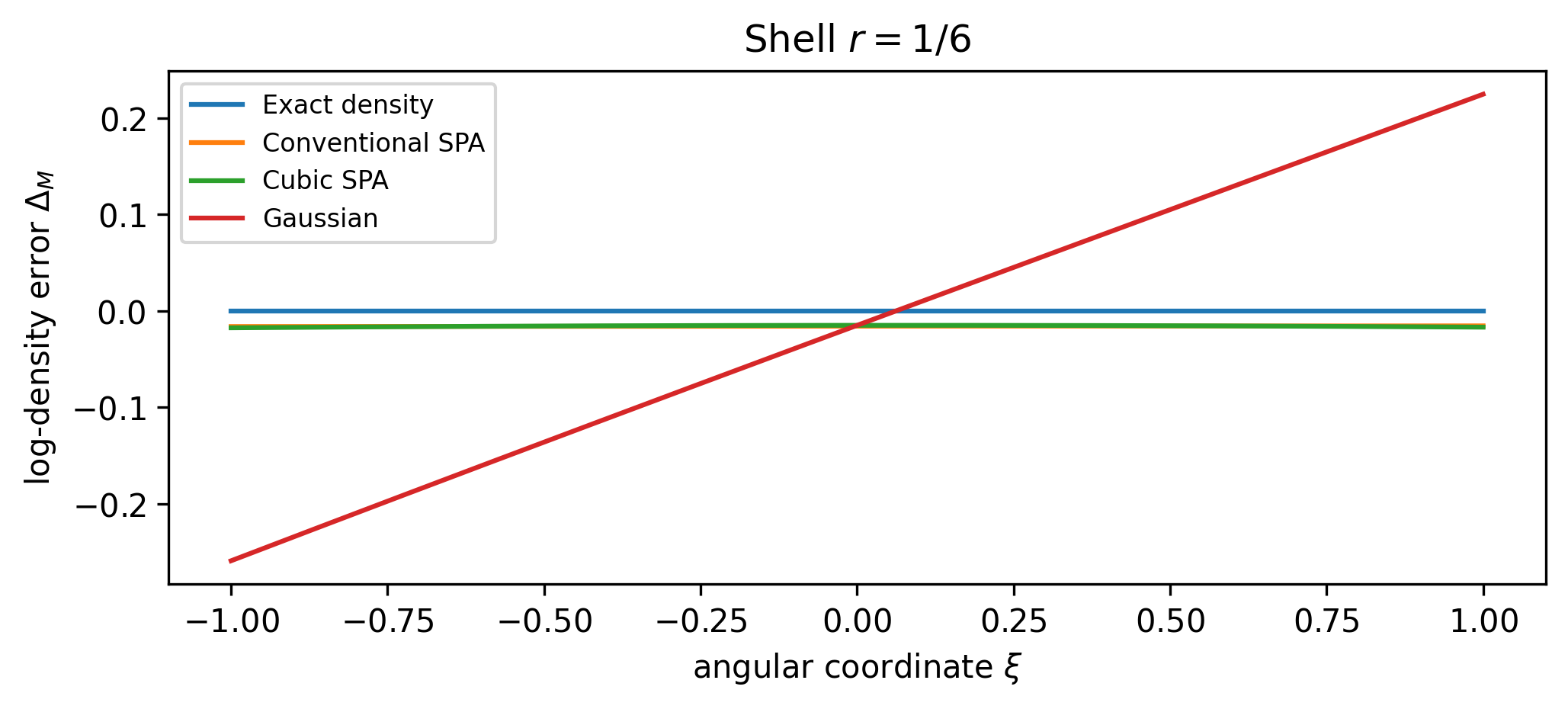}

\vspace{1em}

\includegraphics[width=.75\textwidth]{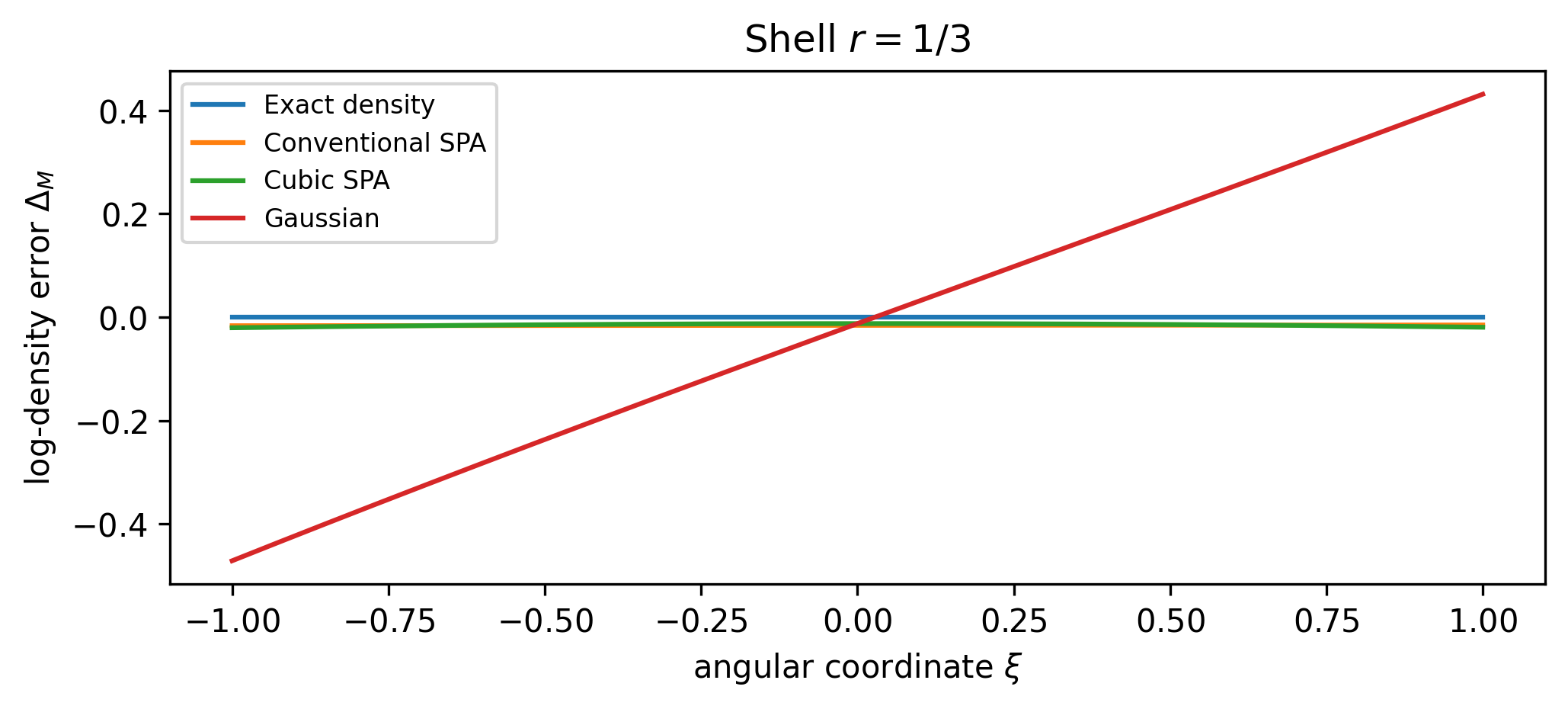}
\caption{Signed log-density errors \eqref{image mixture log errors} of
the conventional SPA, cubic SPA, and Gaussian approximation as functions
of the angular coordinate \eqref{image mixture angular coordinate}.
Top: the shell $\mathcal S_{1/6}$. Bottom: the shell
$\mathcal S_{1/3}$. The exact density corresponds to the horizontal
zero line.}
\label{fig:image mixture shell inner}
\end{figure}

\begin{figure}[h]
\centering
\includegraphics[width=.75\textwidth]{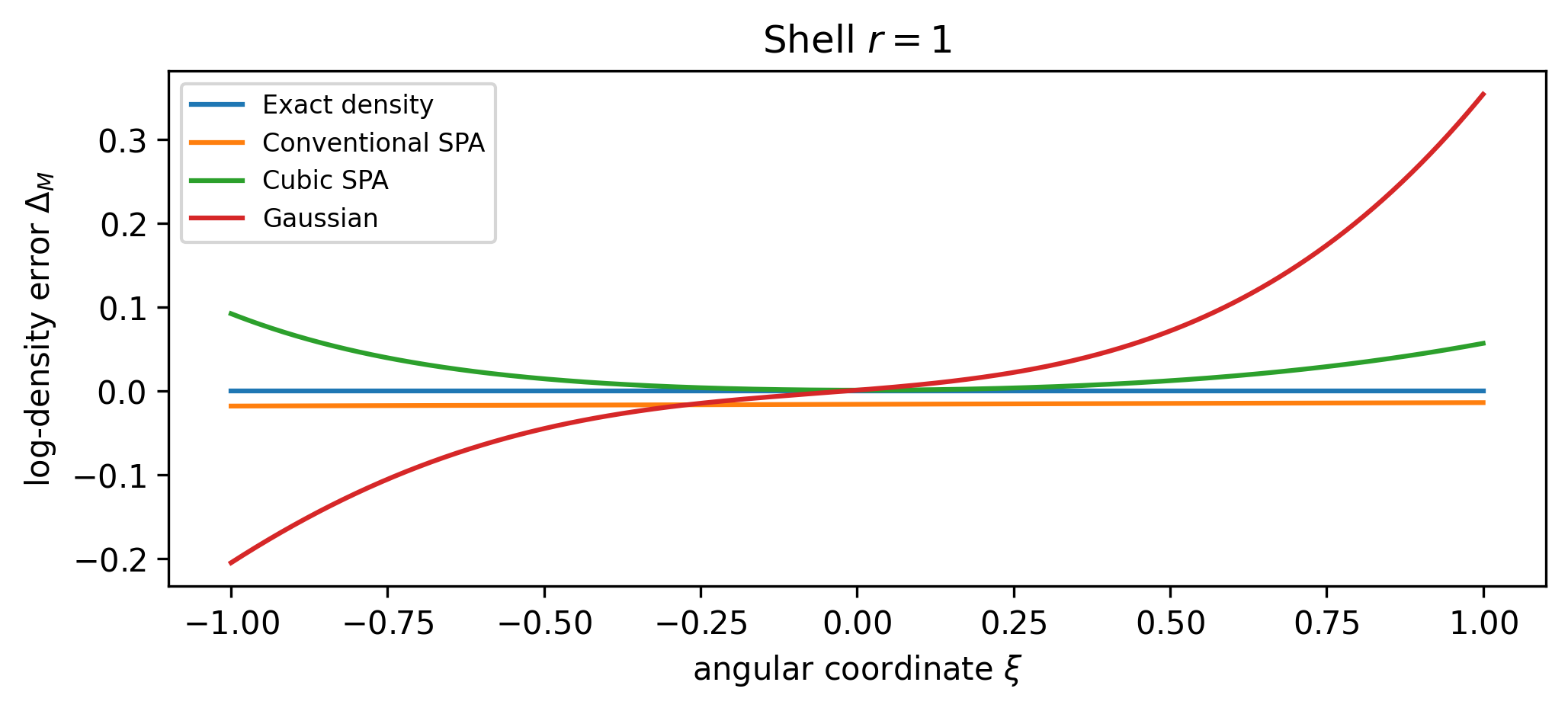}

\vspace{1em}

\includegraphics[width=.75\textwidth]{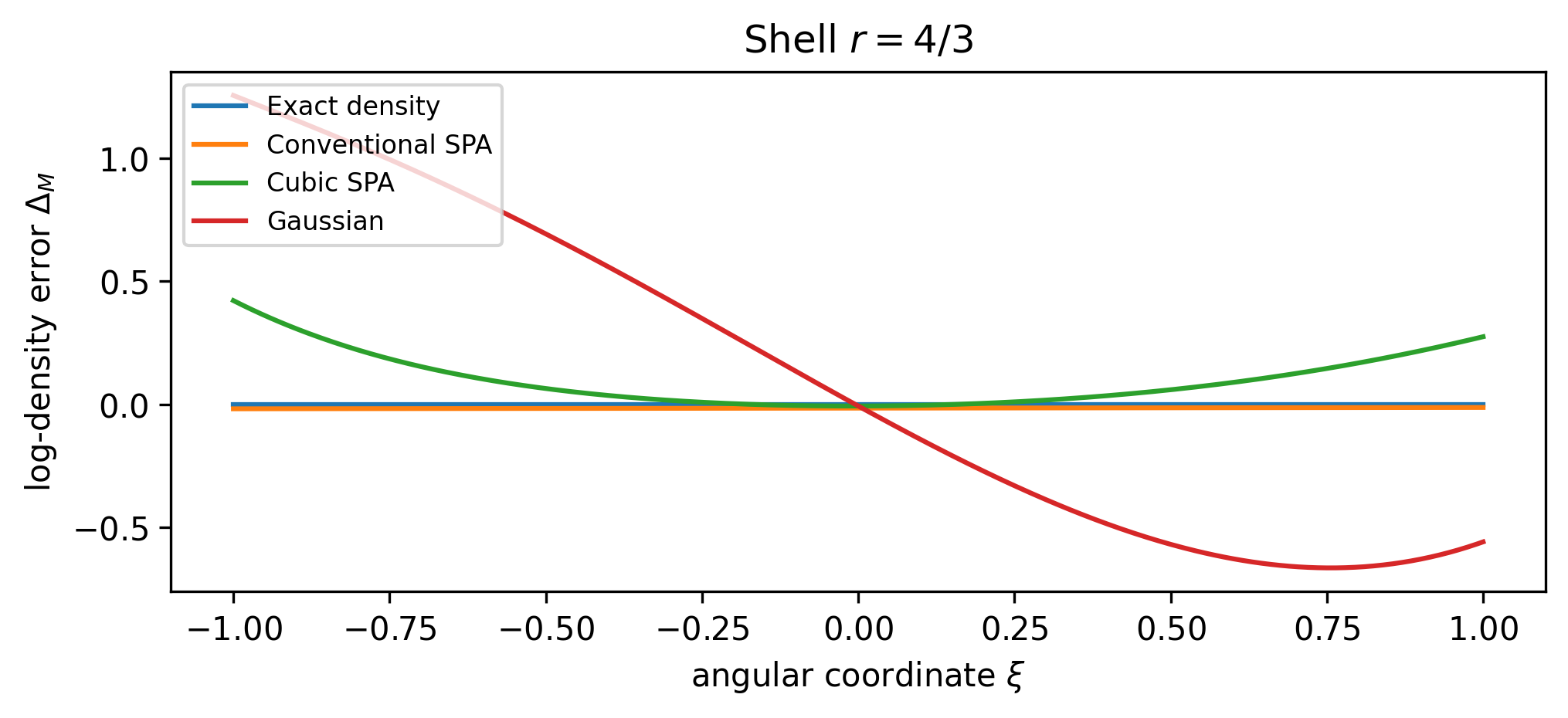}
\caption{Signed log-density errors \eqref{image mixture log errors} of
the conventional SPA, cubic SPA, and Gaussian approximation as functions
of the angular coordinate \eqref{image mixture angular coordinate}.
Top: the shell $\mathcal S_1$. Bottom: the shell
$\mathcal S_{4/3}$. The exact density corresponds to the horizontal
zero line.}
\label{fig:image mixture shell outer}
\end{figure}

Figures~\ref{fig:image mixture shell inner} and \ref{fig:image mixture shell outer} show the $\Delta_M(a)$ profiles as functions of $\xi$ for $r=1/6,1/3,1,$ and $4/3$. The first two are inner central shells, $r=1$ is the natural Gaussian radial scale, and $r=4/3$ probes somewhat farther from the center.

The conventional-SPA error is nearly independent of direction on all four shells and remains close to its value at the central saddlepoint. Indeed, since $\Delta_{\rm S}(a)=-\log I(a)$ and $I(a)$ is close to one, the near constancy of $\Delta_{\rm S}$ is equivalent to the near
constancy of $I(a)-1$. Thus the single Monte Carlo estimate $\widehat J_M(0)-1=1.5935\times10^{-2}$ remains representative of the conventional-SPA error throughout the displayed region, which is qualitatively consistent with the refinement of Theorem~\ref{thm:improved SPA}.

The central expansions exhibit substantially stronger directional dependence. On the two inner shells, the cubic approximation remains close to the conventional SPA, whereas the Gaussian error varies almost linearly with $\xi$. This linear behavior has a simple local explanation. On a fixed shell $S_r$, the Gaussian log density $\log\rho_n^{\rm G}(a)$ is independent of direction, since $\|\Sigma^{-1/2}a\|$ is fixed. Hence the directional dependence of $\Delta_{\rm G}(a)$ comes entirely from the exact log density. By rotational
symmetry in the whitened coordinates, near $a=0$,
\be
\log\rho_n(a)
=
\log\rho_n(0)
+c\langle e,\Sigma^{-1/2}a\rangle
+O(\|\Sigma^{-1/2}a\|^2)
\ee
for some constant $c$. Consequently, on $S_r$ the leading directional term is linear in $\xi$,
with a coefficient proportional to $r(d/n)^{1/2}$; see \eqref{image mixture MC shells} and \eqref{image mixture angular coordinate}. Thus the nearly straight Gaussian profiles on the inner shells and the mode displacement in Table~\ref{tab:image mixture modes} have the same local origin. At $r=1$, the cubic correction still removes most of the directional error of the Gaussian approximation. On the outer shell $r=4/3$, both central expansions deteriorate in directions nearly parallel or antiparallel to $e$, with the Gaussian approximation deteriorating much more strongly. 


We finally illustrate the use of SPA in one particular inverse problem. For the true image, we start from a fixed cross-shaped pattern on the $6\times6$ pixel grid and project it onto the
$4\times4$ DCT model space in \eqref{image mixture DCT}. The resulting image is then normalized so that
\be
\max_x |m_{\theta_0}(x)|=1.5/\sqrt n.
\ee
This scaling keeps the signal amplitude comparable to the residual noise level in the averaged image, which is of order $n^{-1/2}$, and hence makes the reconstruction problem nontrivial. The true image and a single realization of all $1000$ exposures are kept fixed throughout the experiment.

\begin{figure}[t]
\centering
\includegraphics[width=.98\textwidth]{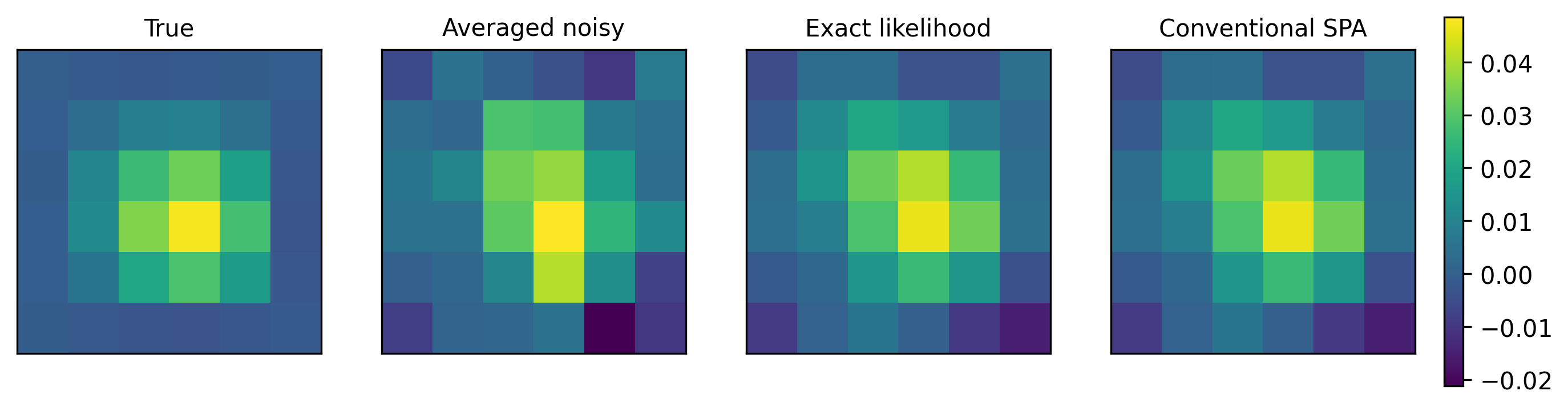}
\caption{True image, averaged noisy image, and reconstructions obtained
by maximizing the exact likelihood and the conventional SPA. The
reconstructions obtained from the cubic and Gaussian approximations are
visually indistinguishable at this scale from the two reconstructions
shown and are therefore omitted.}
\label{fig:image mixture reconstructions}
\end{figure}

Figure~\ref{fig:image mixture reconstructions} shows the true image, the
averaged noisy image, and the reconstructions obtained by maximizing
the exact likelihood and the conventional SPA. The reconstructions
obtained from the cubic and Gaussian approximations are visually
indistinguishable at this scale from the exact-likelihood and
conventional-SPA reconstructions shown in the figure and are therefore
omitted. We report their numerical errors below. The relative $L^2$
errors with respect to $m_{\theta_0}$ are
\be\label{image mixture reconstruction errors}
\frac{\|\bar y-m_{\theta_0}\|_2}{\|m_{\theta_0}\|_2}
=
0.53307,
\quad
\frac{\|\widehat m_M-m_{\theta_0}\|_2}{\|m_{\theta_0}\|_2}
=
0.38251,\quad
0.38252,\quad
0.38251,\quad
0.38549,
\ee
where $M={\rm true},{\rm S},{\rm C},{\rm G}$, respectively, and
\be
\widehat\theta^M:=\arg\max_{\theta\in\br^{16}}L_n^M(\theta;\bar y),
\quad
\widehat m^M:=B\widehat\theta^M.
\ee
Recall that $B$ is defined in \eqref{image mixture DCT}.

We see that all three approximate-likelihood reconstructions give nearly the same point estimate as the exact likelihood in this realization, with the conventional and cubic SPAs agreeing especially closely with it. This indicates that, for this realization, the likelihood maximizer is relatively insensitive to the density-approximation errors. This observation is
consistent with general fixed-dimensional results on the accuracy of saddlepoint maximum-likelihood estimators \cite{goodman2022}.

\appendix

\section{Auxiliary results}

\subsection{Bound on the gradient}
Pick any $f\in C^3(\br^d)$ such that
\be
f(0)=0,\quad \nabla f(0)=0,\quad \nabla^2 f(0)=0.
\ee
Letting $g(\la):=\nabla f(\la x)\cdot u$ for any $u\in S^{d-1}$, we have $g(0)=0$, $g^\prime(0)=\nabla^2 f(0)[x,u]=0$, and $g^{\prime\prime}(\la)=\nabla^3 f(\la x)[x,x,u]$. Taylor’s formula with integral remainder yields
\be\label{g1}
g(1)=\int_0^1 (1-\la)\,g^{\prime\prime}(\la)\,\dd\la =\int_0^1 (1-\la)\, \nabla^3 f(\la x)[x,x,u]\,\dd\la.
\ee
Taking the supremum over $u\in S^{d-1}$ yields
\be\label{grad bnd}
\Vert\nabla f(x)\Vert
\le \tfrac12 \sup_{\la\in[0,1]}\Vert \nabla^3 f(\la x)\Vert_{\op}\Vert x\Vert^2.
\ee

\subsection{A property of the Legendre transform}\label{ssec:legendre}

Consider the equation $\nabla\varphi(\tau)=a$. In the setting of Subsection~\ref{ssec:local CLT}, $X$ is isotropic, so $\nabla\varphi(0)=0$ and $\nabla^2\varphi(0)=I_d$. Similarly to \eqref{g1}, 
\be\label{gradphi}
\langle\nabla \varphi(\tau),u\rangle=\langle\tau,u\rangle+\int_0^1 (1-\la)\, \nabla^3 \varphi(\la \tau)[\tau,\tau,u]\,\dd\la=\langle a,u\rangle.
\ee
This implies
\be\label{atau dist}
\Vert \tau-a\Vert\le \tfrac12\max_{\la\in[0,1]}\Vert\nabla^3 \varphi(\la \tau)\Vert_{\op} \Vert \tau\Vert^2.
\ee
Assuming $\max_{\la\in[0,1]}\Vert\nabla^3 \varphi(\la \tau)\Vert_{\op} \Vert \tau\Vert\le1$, we get $\Vert \tau\Vert\le 2\Vert a\Vert$. 

Further, a simple manipulation shows that
\be\label{gradphi err}
\varphi^*(a)-\tfrac12 \Vert a\Vert^2=-\tfrac12 \Vert \tau-a\Vert^2-
\tfrac16\nabla^3 \varphi(\la \tau)[\tau,\tau,\tau]
\ee
for some $\la\in(0,1)$. Combining with the preceding results proves \eqref{poly-bound}. 

Note that \eqref{C3 def}  implies that the assumption $\max_{\la\in[0,1]}\Vert\nabla^3 \varphi(\la \tau)\Vert_{\op} \Vert \tau\Vert\le1$ is satisfied. This and the existence of a solution $\tau$ in the set $\Vert \tau\Vert\le 2\Vert a\Vert$ for any $a$ that satisfies \eqref{C3 def} follows by removing $u$ from \eqref{gradphi} and applying Brouwer's fixed-point theorem to the resulting equation. 

Indeed, define the matrix function:
\be\label{matr fn}
B(\tau)=\int_0^1 (1-\la)\, \nabla^3 \varphi(\la \tau)\tau\,\dd\la.
\ee
By construction, $\|B(\tau)\| \le 1/2$ whenever $\|\tau\| \le 2\|a\|$, $a \in U_d$, and \eqref{C3 def} holds. Under these conditions, \eqref{gradphi} is equivalent to the fixed-point equation $\tau = D(\tau)$, where $D(\tau) := (\mathrm{I} + B(\tau))^{-1} a$, defined on the ball $\|\tau\| \le 2\|a\|$. Since $\|(\mathrm{I} + B(\tau))^{-1}\| \le 2$, the map $D(\tau)$ sends this ball into itself. By Brouwer’s fixed-point theorem, a solution exists. Uniqueness follows from the strict positivity of $\nabla^2 \varphi(\tau)$, $\tau\in V_d$, cf. Assumption~\ref{ass:mainI}\eqref{posdef}.

\subsection{Proof of Corollary~\ref{cor:central density expansion}}\label{sec:prf cor expansion}

Set $a=x/\sqrt n$. On the region $\|x\|\le L\sqrt d$,
\be\label{a norm}
\|a\|\le L(d/n)^{1/2}\to0
\ee
when $d^2/n\to0$. 

By Lemma~\ref{lem:aux ass}, there exists a fixed $0<r'\le r$ such that parts \eqref{exist}, \eqref{posdef}, and \eqref{phase} of Assumption~\ref{ass:mainI} hold on $V_d'=\{\|\tau\|<r'\}$, and $c_3,c_4$ are uniformly bounded. Part \eqref{global} holds there by restriction from $\{\|\tau\|<r\}$.

Since \eqref{a norm} and the uniform local bounds on $H(\tau)$ in \eqref{spectral bounds} imply $\|\tau\| \les\|a\|$, the corresponding saddlepoints $\tau(a)$ converge to zero uniformly for $\|x\|\le L\sqrt d$. We used here that $X$ is isotropic and that \eqref{uniform directional moments} holds; together these imply \eqref{spectral bounds}, by the argument following \eqref{spectral bounds}. Hence all $\tau(a)$ belong to $V_d'$ for $d,n$ sufficiently large, and Theorem~\ref{thm:mainI} may be applied on $V_d'$.

Taylor expansion gives
\be\label{central cgf expansion}
\varphi(\tau)
=\tfrac12\|\tau\|^2 +\tfrac16T_3[\tau,\tau,\tau] +O(\|\tau\|^4)
\ee
uniformly in $d$, and therefore
\be\label{central grad expansion}
\nabla\varphi(\tau)
= \tau+\tfrac12T_3[\tau,\tau,\cdot] +O(\|\tau\|^3).
\ee
Solving $\nabla\varphi(\tau)=a$ by one substitution yields
\be\label{central saddle expansion}
\tau=a-\tfrac12T_3[a,a,\cdot]+O(\|a\|^3).
\ee
Substitution into
$\varphi^*(a)=\langle\tau,a\rangle-\varphi(\tau)$ gives
\be\label{central Legendre expansion}
\varphi^*(a)
=\tfrac12\|a\|^2-\tfrac16T_3[a,a,a]+O(\|a\|^4).
\ee
Since $n\|a\|^4=O_L(d^2/n)$, we obtain
\be\label{central Legendre scaled}
n\big(\varphi^*(a)-\tfrac12\|a\|^2\big)
=-\tfrac1{6\sqrt n}T_3[x,x,x]+O_L\big(\tfrac{d^2}{n}\big).
\ee

Similarly,
\be\label{central Hess expansion}
H(\tau)
=
I_d+T_3[\tau,\cdot,\cdot]+O(\|\tau\|^2),
\ee
and \eqref{central saddle expansion} gives
\be
H(\tau)
=
I_d+T_3[a,\cdot,\cdot]+O(\|a\|^2).
\ee
Expanding the logarithm of the determinant around $I_d$ therefore gives
\be\label{central logdet expansion}
\log\det H(\tau)
=\tr T_3[a,\cdot,\cdot]+O(d\|a\|^2)
=\tfrac1{\sqrt n}\tr T_3[x,\cdot,\cdot]+O_L\left(\tfrac{d^2}{n}\right).
\ee

The density of $S_n$ is related to the density $\rho_n$ of the sample mean by
\be
f_n(x)=n^{-d/2}\rho_n(x/\sqrt n).
\ee
Using \eqref{rho a}, \eqref{central Legendre scaled}, \eqref{central logdet expansion}, and Theorem~\ref{thm:mainI} now gives \eqref{central density expansion}.

\subsection{A lower bound for $\|(nH_r)^{-1/2}\|_\infty$}
\label{ssec:TR infinity lower bound}

We show that the factor $d^{1/2}$ in \eqref{TR infinity bound} cannot in general be removed under \eqref{VG param bnds}. It suffices to consider $z=\tau=0$ and
$\Sigma=I_d$. Since $A(0)=1$, \eqref{VG grad Hess} gives
\be
H(0)=s(I_d+bb^\top).
\ee
Choose
\be
b=2\big(1,(d-1)^{-1/2},\ldots,(d-1)^{-1/2}\big)^\top.
\ee
Then $\|b\|^2=8$, so this choice satisfies \eqref{VG param bnds} uniformly in $d$. By an easy calculation,
\be
(I_d+bb^\top)^{-1/2}
=
I_d-\tfrac1{12}bb^\top.
\ee
Since $H_r(0)=H(0)$, it follows that
\be\label{TR infinity exact}
(nH_r(0))^{-1/2}
= (ns)^{-1/2} \big(I_d-\tfrac1{12}bb^\top\big).
\ee

Consider the first row of the matrix in \eqref{TR infinity exact}.
Its diagonal entry equals
\be
1-\tfrac1{12}b_1^2=\tfrac23,
\ee
while each of its $d-1$ off-diagonal entries has absolute value
\be
\tfrac1{12}|b_1b_j|
= \tfrac{1}{3\sqrt{d-1}}, \quad j=2,\ldots,d.
\ee
Therefore
\be
\|(nH_r(0))^{-1/2}\|_\infty
\ge
\tfrac1{\sqrt{ns}}
\left(
\tfrac23+\tfrac13\sqrt{d-1}
\right)
\gtrsim
(d/n)^{1/2}.
\ee


\section{Derivation of the cubic likelihood in Section~\ref{sec:image application}}
\label{app:image cubic derivation}

We give the details of the calculations leading to \eqref{image mixture cubic}, \eqref{q and T3}. From \eqref{image mixture noise},
\be\label{app image X}
X=\sigma\big(S\mu+\sqrt{c_S}G_\al Z\big),\quad Z\sim\CN(0,I_d).
\ee
Here $S$ has the distribution in \eqref{image mixture S}, 
\be
\Eb S=0,\quad \Eb S^2=1,\quad \Eb c_S=1,
\ee
and $Z$ is independent of $S$. It follows immediately that $\Eb X=0$ and
\be\label{VarX}
\Var(X)=\sigma^2(C_\al+\mu\mu^\top)=\Sigma,
\ee
as in \eqref{image mixture Sigma}.

Since $X$ is centered,
\be
T_3=\nabla^3\varphi(0)=\Eb X^{\otimes3}.
\ee
Expanding $X^{\otimes3}$ using \eqref{app image X}, all terms containing an odd number of Gaussian factors vanish. Therefore,
\bs\label{app image T3 derivation}
T_3
&=
\sigma^3\left[
(\Eb S^3)\mu^{\otimes3}
+
\Eb(Sc_S)\operatorname{Sym}(\mu\otimes C_\al)
\right]\\
&=
\sigma^3\left[
\kappa_3\mu^{\otimes3}
+
\beta\operatorname{Sym}(\mu\otimes C_\al)
\right],
\es
where 
\bs
\operatorname{Sym}(\mu\otimes C_\al)[u,u,u]
=&
3\langle\mu,u\rangle\langle u,C_\al u\rangle,\quad u\in\br^d,\\
\operatorname{Sym}(\mu\otimes C)_{ijk}
=&\mu_iC_{jk}+\mu_jC_{ik}+\mu_kC_{ij}.
\es
This proves the first equation in \eqref{q and T3}.

We next derive the vector $q$ in \eqref{q and T3}. This vector appears in the expansion of the determinant term below and it is defined by the property
\be\label{tr id}
\operatorname{tr}\{A T_3[u,\cdot,\cdot]\}= 
\sum_{i=1}^d\Big(\sum_{j,k=1}^d A_{jk}(T_3)_{ijk}\Big) u_i
=\langle q,u\rangle,
\ee
where $A=\Sigma^{-1}$. We have
\be
\sum_{j,k=1}^d \mu_i\mu_j\mu_k A_{jk}
= \langle\mu,A\mu\rangle\mu_i
\ee
and
\bs\label{sym A}
&\sum_{j,k=1}^d
\big(
\mu_i(C_\al)_{jk}
+
\mu_j(C_\al)_{ik}
+
\mu_k(C_\al)_{ij}
\big)A_{jk}\\
&\hspace{1cm}
=
\tr(C_\al A)\mu_i
+
2(C_\al A\mu)_i.
\es
Combining \eqref{app image T3 derivation} with \eqref{tr id}--\eqref{sym A} therefore gives
\be
q
=
\sigma^3
\left[
\kappa_3\langle\mu,A\mu\rangle\mu
+
\beta\left\{
\tr(C_\al A)\mu+2C_\al A\mu
\right\}
\right],
\ee
which is \eqref{q and T3}.

It remains to derive the cubic approximation to the likelihood. This is the
non-isotropic version of the calculation leading to 
\eqref{central saddle expansion}--\eqref{central logdet expansion}. Indeed, since $\Eb X=0$ and $\nabla^2\varphi(0)=\Sigma$, similar Taylor expansions, with $A=\Sigma^{-1}$ and $u=Aa$, give 
\bs\label{app image general expansions}
\tau
&=
u-\tfrac12 A T_3[u,u,\cdot]+O(\|a\|^3),\\
\varphi^*(a)
&=
\tfrac12\langle a,Aa\rangle
-\tfrac16T_3[u,u,u]+O(\|a\|^4).
\es
Similarly,
\be
H(\tau)
=
\Sigma+T_3[u,\cdot,\cdot]+O(\|a\|^2),
\ee
and hence
\bs\label{app image determinant expansion}
\log\det H(\tau)
&=
\log\det\Sigma
+
\tr\{A\,T_3[u,\cdot,\cdot]\}
+
O(d\|a\|^2)\\
&=
\log\det\Sigma+\langle q,u\rangle+O(d\|a\|^2),
\es
where the last equality follows from \eqref{tr id}.

The logarithm of the conventional SPA for the density of $n^{-1}\sum_{j=1}^nX_j$ is
\be
\tfrac d2\log n-\tfrac d2\log(2\pi)
-\tfrac12\log\det H(\tau)-n\varphi^*(a).
\ee
Substitution of the preceding expansions for $\varphi^*(a)$ and $H(\tau)$ gives \be\label{app image cubic likelihood derivation}
C_{n,\Sigma}
-\tfrac n2\langle a,Aa\rangle
+\tfrac n6T_3[u,u,u]
-\tfrac12\langle q,u\rangle,
\ee
which is \eqref{image mixture cubic}.

Finally, in the central region \eqref{image mixture central region},
$\|a\|=O_R((d/n)^{1/2})$ under the uniform spectral bounds on $\Sigma$.
Consequently both $nO(\|a\|^4)$ and $O(d\|a\|^2)$ are
$O_R(d^2/n)$. Together with Theorem~\ref{thm:mainI}, this gives
\eqref{image mixture cubic accuracy}.

\section{Evaluation of $c_3$ and $c_4$ in the numerical experiment}
\label{app:image c3c4}

We explain the numerical evaluation of the derivative constants used to assess the bound in Theorem~\ref{thm:mainI} at $a=0$. The corresponding saddlepoint is $\tau=0$. Set
\be
M:=H(0)^{-1/2},\quad v:=M\mu.
\ee
Replacing $\tau$ with $iMt$ in \eqref{image mixture cgf} gives
\be\label{app image cgf Mt}
\varphi(iMt) = F\left(\langle v,t\rangle,\langle Mt,C_\al Mt\rangle\right)
\ee
for an explicit scalar function $F$. We have $H(0)=\sigma^2(C_\al+\mu\mu^\top)$; see \eqref{image mixture Sigma} and, for its derivation, \eqref{VarX}. Therefore, 
\be
\sigma^2\big(MC_\al M+vv^\top\big)=I_d.
\ee
Consequently,
\be\label{app image cgf reduced}
\varphi(iMt)
=
F\left(
\langle v,t\rangle,
\sigma^{-2}\|t\|^2-\langle v,t\rangle^2
\right).
\ee

To compute the tensor operator norms at a fixed $t$, let $\Theta\in S^{d-1}$ and replace $t$ with $t+h\Theta$ in \eqref{app image cgf reduced}. Set
\be
x:=\langle v,t\rangle,\quad r:=\|t\|,\quad
u:=\langle v,\Theta\rangle,\quad w:=\langle t,\Theta\rangle.
\ee
Then
\be\label{app image directional F}
\varphi(iM(t+h\Theta))
=
F\left(
x+hu,\,
\sigma^{-2}(r^2+2hw+h^2)-(x+hu)^2
\right).
\ee
Thus the directional derivatives of orders three and four with respect to $h$ at $h=0$ depend only on the four scalar quantities $x,r,u,w$. For each fixed $t$, the computation of the tensor operator norm therefore reduces to a two-dimensional optimization over the geometrically feasible pairs $(u,w)$, with $|u|\le \|v\|$ and $|w|\le r$, while the supremum over the ball $\{\|t\|<2.5\e^{1/2}\}$ in \eqref{c3c4} reduces to a two-dimensional optimization in $(x,r)$, with $|x|\le \|v\|r$. Hence the full computation of $c_{3,0}$ and $c_{4,0}$ is a four-dimensional compact optimization problem.

\section*{Declaration of generative AI use}

OpenAI's ChatGPT was used in the preparation of this manuscript for substantial
assistance with the Python implementation of the numerical experiments, for routine
algebraic calculations, and for English-language editing and stylistic suggestions.
The Python code was reviewed, tested, and validated by the author. The mathematical
arguments and conclusions are the author's own, and the author takes full responsibility
for the content of the manuscript and the reported numerical results.

\bibliographystyle{abbrv}
\bibliography{My_Collection}

@article{anderson1998,
  title = {Edgeworth Expansions in Very-High-Dimensional Problems},
  author = {Anderson, Niall H. and Hall, Peter and Titterington, D.M.},
  year = 1998,
  month = jul,
  journal = {Journal of Statistical Planning and Inference},
  volume = {70},
  number = {1},
  pages = {1--18},
  issn = {03783758},
  doi = {10.1016/S0378-3758(97)00171-7},
  urldate = {2026-09-23},
  copyright = {https://www.elsevier.com/tdm/userlicense/1.0/},
  langid = {english}
}

@book{bakry2014,
  title = {Analysis and {{Geometry}} of {{Markov Diffusion Operators}}},
  author = {Bakry, Dominique and Gentil, Ivan and Ledoux, Michel},
  year = 2014,
  series = {Grundlehren Der Mathematischen {{Wissenschaften}}},
  volume = {348},
  publisher = {Springer International Publishing},
  address = {Cham},
  doi = {10.1007/978-3-319-00227-9},
  urldate = {2025-10-13},
  copyright = {https://www.springernature.com/gp/researchers/text-and-data-mining},
  isbn = {978-3-319-00226-2 978-3-319-00227-9},
  langid = {english}
}

@article{belloni2014,
  title = {Posterior Inference in Curved Exponential Families under Increasing Dimensions},
  author = {Belloni, Alexandre and Chernozhukov, Victor},
  year = 2014,
  month = jun,
  journal = {The Econometrics Journal},
  volume = {17},
  number = {2},
  pages = {S75-S100},
  issn = {1368-4221, 1368-423X},
  doi = {10.1111/ectj.12027},
  urldate = {2026-09-18},
  copyright = {http://creativecommons.org/licenses/by/3.0/},
  langid = {english}
}

@article{bentkus1996,
  title = {Optimal Rates of Convergence in the {{CLT}} for Quadratic Forms},
  author = {Bentkus, V. and G{\"o}tze, F.},
  year = 1996,
  month = jan,
  journal = {The Annals of Probability},
  volume = {24},
  number = {1},
  issn = {0091-1798},
  doi = {10.1214/aop/1042644727},
  urldate = {2026-09-10}
}

@article{bentkus2003,
  title = {On the Dependence of the {{Berry}}--{{Esseen}} Bound on Dimension},
  author = {Bentkus, V.},
  year = 2003,
  month = may,
  journal = {Journal of Statistical Planning and Inference},
  volume = {113},
  number = {2},
  pages = {385--402},
  issn = {03783758},
  doi = {10.1016/S0378-3758(02)00094-0},
  urldate = {2026-09-10},
  copyright = {https://www.elsevier.com/tdm/userlicense/1.0/},
  langid = {english}
}

@article{bentkus2005,
  title = {A {{Lyapunov-type Bound}} in {\emph{ }}{{{\emph{R}}}}{\emph{{\textsuperscript{d}}}}{\emph{ }}},
  author = {Bentkus, V.},
  year = 2005,
  month = jan,
  journal = {Theory of Probability \& Its Applications},
  volume = {49},
  number = {2},
  pages = {311--323},
  issn = {0040-585X, 1095-7219},
  doi = {10.1137/S0040585X97981123},
  urldate = {2026-09-10},
  langid = {english}
}

@book{bhattacharya2010,
  title = {Normal Approximation and Asymptotic Expansions},
  author = {Bhattacharya, R. N. and Ranga Rao, R.},
  year = 2010,
  series = {Classics in Applied Mathematics},
  number = {64},
  publisher = {{Society for Industrial and Applied Mathematics}},
  address = {Philadelphia},
  isbn = {978-0-89871-897-3},
  lccn = {QA273.67 .B48 2010}
}

@book{billingsley2012,
  title = {Probability and Measure},
  author = {Billingsley, Patrick},
  year = 2012,
  series = {Wiley Series in Probability and Statistics},
  edition = {Anniversary},
  publisher = {Wiley},
  address = {Hoboken, New Jersey},
  isbn = {978-1-118-12237-2},
  langid = {english}
}

@article{bobkov2000,
  title = {From {{Brunn-Minkowski}} to {{Brascamp-Lieb}} and to Logarithmic {{Sobolev}} Inequalities},
  author = {Bobkov, S.G. and Ledoux, M.},
  year = 2000,
  month = dec,
  journal = {Geometric and Functional Analysis},
  volume = {10},
  number = {5},
  pages = {1028--1052},
  issn = {1016-443X, 1420-8970},
  doi = {10.1007/PL00001645},
  urldate = {2026-09-11},
  langid = {english}
}

@article{bobkov2012,
  title = {Bounds for Characteristic Functions in Terms of Quantiles and Entropy},
  author = {Bobkov, Sergey and Chistyakov, Gennadiy and G{\"o}tze, Friedrich},
  year = 2012,
  month = jan,
  journal = {Electronic Communications in Probability},
  volume = {17},
  issn = {1083-589X},
  doi = {10.1214/ECP.v17-2053},
  urldate = {2026-09-02}
}

@article{bobkov2025,
  title = {Berry--{{Esseen}} Bounds in the Local Limit Theorems},
  author = {Bobkov, Sergey G. and G{\"o}tze, Friedrich},
  year = 2025,
  month = jan,
  journal = {Lithuanian Mathematical Journal},
  volume = {65},
  number = {1},
  pages = {50--66},
  issn = {0363-1672, 1573-8825},
  doi = {10.1007/s10986-025-09659-1},
  urldate = {2026-09-17},
  langid = {english}
}

@incollection{broda2012,
  title = {Saddlepoint {{Approximations}}: {{A Review}} and {{Some New Applications}}},
  shorttitle = {Saddlepoint {{Approximations}}},
  booktitle = {Handbook of {{Computational Statistics}}},
  author = {Broda, Simon A. and Paolella, Marc S.},
  editor = {Gentle, James E. and H{\"a}rdle, Wolfgang Karl and Mori, Yuichi},
  year = 2012,
  pages = {953--983},
  publisher = {Springer Berlin Heidelberg},
  address = {Berlin, Heidelberg},
  doi = {10.1007/978-3-642-21551-3_32},
  urldate = {2025-10-16},
  isbn = {978-3-642-21550-6 978-3-642-21551-3},
  langid = {english}
}

@book{butler2007,
  title = {Saddlepoint Approximations with Applications},
  author = {Butler, Ronald W.},
  year = 2007,
  series = {Cambridge Series in Statistical and Probabilistic Mathematics},
  publisher = {Cambridge University Press},
  address = {Cambridge},
  isbn = {978-0-521-87250-8 978-0-511-34074-1},
  langid = {english}
}

@incollection{chakraborty2025,
  title = {The {{Central Limit Theorem}} in {{High-Dimensions}} and {{Its Applications}}},
  booktitle = {Frontiers of {{Statistics}} and {{Data Science}}},
  author = {Chakraborty, Nilanjan and Lahiri, S. N.},
  editor = {Ghosal, Subhashis and Roy, Anindya},
  year = 2025,
  pages = {81--100},
  publisher = {Springer Nature Singapore},
  address = {Singapore},
  doi = {10.1007/978-981-96-0742-6_4},
  urldate = {2025-10-15},
  isbn = {978-981-96-0741-9 978-981-96-0742-6},
  langid = {english}
}

@article{chernozhukov2017,
  title = {Central Limit Theorems and Bootstrap in High Dimensions},
  author = {Chernozhukov, Victor and Chetverikov, Denis and Kato, Kengo},
  year = 2017,
  month = jul,
  journal = {The Annals of Probability},
  volume = {45},
  number = {4},
  issn = {0091-1798},
  doi = {10.1214/16-AOP1113},
  urldate = {2026-09-10}
}

@article{Daniels1954,
  title = {Saddlepoint {{Approximations}} in {{Statistics}}},
  author = {Daniels, H. E.},
  year = 1954,
  journal = {The Annals of Mathematical Statistics},
  volume = {25},
  number = {4},
  pages = {631--650},
  issn = {0003-4851},
  doi = {10.1214/aoms/1177728652}
}

@article{goodman2022,
  title = {Asymptotic Accuracy of the Saddlepoint Approximation for Maximum Likelihood Estimation},
  author = {Goodman, Jesse},
  year = 2022,
  month = aug,
  journal = {The Annals of Statistics},
  volume = {50},
  number = {4},
  issn = {0090-5364},
  doi = {10.1214/22-AOS2169},
  urldate = {2026-09-29}
}

@article{gotze1991,
  title = {On the {{Rate}} of {{Convergence}} in the {{Multivariate CLT}}},
  author = {G{\"o}tze, F.},
  year = 1991,
  month = apr,
  journal = {The Annals of Probability},
  volume = {19},
  number = {2},
  issn = {0091-1798},
  doi = {10.1214/aop/1176990448},
  urldate = {2026-09-10}
}

@article{haroske2009,
  title = {Besov Spaces with Positive Smoothness on $\mathbb{R}^n$, Embeddings and Growth Envelopes},
  author = {Haroske, Dorothee D. and Schneider, Cornelia},
  year = 2009,
  month = dec,
  journal = {Journal of Approximation Theory},
  volume = {161},
  number = {2},
  pages = {723--747},
  issn = {00219045},
  doi = {10.1016/j.jat.2008.12.004},
  urldate = {2026-09-02},
  langid = {english}
}

@article{jensen1991,
  title = {Uniform {{Saddlepoint Approximations}} and {{Log-Concave Densities}}},
  author = {Jensen, Jens Ledet},
  year = 1991,
  month = sep,
  journal = {Journal of the Royal Statistical Society Series B: Statistical Methodology},
  volume = {53},
  number = {1},
  pages = {157--172},
  issn = {1369-7412, 1467-9868},
  doi = {10.1111/j.2517-6161.1991.tb01815.x},
  urldate = {2026-09-22},
  copyright = {https://academic.oup.com/journals/pages/open\_access/funder\_policies/chorus/standard\_publication\_model},
  langid = {english}
}

@book{jensen1995,
  title = {Saddlepoint {{Approximations}}},
  author = {Jensen, Jens Ledet},
  year = 1995,
  month = may,
  publisher = {Oxford University Press},
  address = {Oxford},
  doi = {10.1093/oso/9780198522959.001.0001},
  urldate = {2026-09-22},
  isbn = {978-0-19-852295-9 978-1-383-02384-8},
  langid = {english}
}

@article{jensen2021,
  title = {On the {{Use}} of {{Saddlepoint Approximations}} in {{High Dimensional Inference}}},
  author = {Jensen, Jens Ledet},
  year = 2021,
  month = feb,
  journal = {Sankhya A},
  volume = {83},
  number = {1},
  pages = {379--392},
  issn = {0976-836X, 0976-8378},
  doi = {10.1007/s13171-019-00188-x},
  urldate = {2026-09-24},
  langid = {english}
}

@article{kasprzak2025,
  title = {How Good Is Your {{Laplace}} Approximation of the {{Bayesian}} Posterior? {{Finite-sample}} Computable Error Bounds for a Variety of Useful Divergences},
  author = {Kasprzak, Mikolaj J and Giordano, Ryan and Broderick, Tamara},
  year = 2025,
  journal = {Journal of Machine Learning Research},
  volume = {26},
  pages = {1--81},
  langid = {english}
}

@misc{katsevich2024b,
  title = {The {{Laplace}} Approximation Accuracy in High Dimensions: A Refined Analysis and New Skew Adjustment},
  shorttitle = {The {{Laplace}} Approximation Accuracy in High Dimensions},
  author = {Katsevich, Anya},
  year = 2024,
  month = jun,
  number = {arXiv:2306.07262},
  eprint = {2306.07262},
  publisher = {arXiv},
  doi = {10.48550/arXiv.2306.07262},
  urldate = {2025-10-15},
  archiveprefix = {arXiv},
  langid = {english},
  note = {arXiv:2306.07262 [math]}
}

@misc{katsevich2025a,
  title = {The {{Laplace}} Asymptotic Expansion in High Dimensions},
  author = {Katsevich, Anya},
  year = 2025,
  month = jun,
  number = {arXiv:2406.12706},
  eprint = {2406.12706},
  publisher = {arXiv},
  doi = {10.48550/arXiv.2406.12706},
  urldate = {2025-10-15},
  archiveprefix = {arXiv},
  langid = {english},
  note = {arXiv:2406.12706 [math]}
}

@article{katsevich2025c,
  title = {Improved Dimension Dependence in the {{Bernstein}}--von {{Mises}} Theorem via a New {{Laplace}} Approximation Bound},
  author = {Katsevich, Anya},
  year = 2025,
  journal = {Information and Inference: A Journal of the IMA},
  volume = {14},
  number = {3},
  pages = {iaaf020},
  issn = {2049-8772},
  doi = {10.1093/imaiai/iaaf020},
  urldate = {2026-07-13},
  copyright = {https://academic.oup.com/pages/standard-publication-reuse-rights},
  langid = {english}
}

@misc{katsevich2026a,
  title = {High-Dimensional {{Laplace}} Asymptotics up to the Concentration Threshold},
  author = {Katsevich, Alexander and Katsevich, Anya},
  year = 2026,
  publisher = {arXiv},
  doi = {10.48550/ARXIV.2602.23151},
  urldate = {2026-07-01},
  copyright = {arXiv.org perpetual, non-exclusive license},
  note = {arXiv:2602.23151v2}
}

@article{Kolassa2003,
  title = {Multivariate Saddlepoint Tail Probability Approximations},
  author = {Kolassa, John E.},
  year = 2003,
  journal = {Annals of Statistics},
  volume = {31},
  number = {1},
  pages = {274--286},
  issn = {00905364},
  doi = {10.1214/aos/1046294465}
}

@article{Kolassa2010,
  title = {Multivariate Saddlepoint Approximations in Tail Probability and Conditional Inference},
  author = {Kolassa, John and Li, Jixin},
  year = 2010,
  journal = {Bernoulli},
  volume = {16},
  number = {4},
  eprint = {1011.5775v1},
  pages = {1191--1207},
  issn = {13507265},
  doi = {10.3150/09-BEJ237},
  archiveprefix = {arXiv}
}

@article{kozubowski2013,
  title = {Multivariate Generalized {{Laplace}} Distribution and Related Random Fields},
  author = {Kozubowski, Tomasz J. and Podg{\'o}rski, Krzysztof and Rychlik, Igor},
  year = 2013,
  month = jan,
  journal = {Journal of Multivariate Analysis},
  volume = {113},
  pages = {59--72},
  issn = {0047259X},
  doi = {10.1016/j.jmva.2012.02.010},
  urldate = {2026-09-07},
  copyright = {https://www.elsevier.com/tdm/userlicense/1.0/},
  langid = {english}
}

@article{lewis2009,
  title = {Galaxy Shear Estimation from Stacked Images},
  author = {Lewis, Antony},
  year = 2009,
  month = sep,
  journal = {Monthly Notices of the Royal Astronomical Society},
  volume = {398},
  number = {1},
  pages = {471--476},
  issn = {00358711, 13652966},
  doi = {10.1111/j.1365-2966.2009.15161.x},
  urldate = {2026-09-07},
  langid = {english}
}

@article{LV2007,
  title = {The Geometry of Logconcave Functions and Sampling Algorithms},
  author = {Lov{\'a}sz, L{\'a}szl{\'o} and Vempala, Santosh},
  year = 2007,
  journal = {Random Structures and Algorithms},
  volume = {30},
  number = {3},
  pages = {307--358},
  issn = {10429832},
  doi = {10.1002/rsa.20135}
}

@inproceedings{mauer2009,
  title = {Empirical {{Bernstein}} Bounds and Sample Variance Penalization},
  booktitle = {22nd {{Annual Conference}} on {{Learning Theory}} ({{COLT}})},
  author = {Maurer, A and Pontil, M},
  year = 2009
}

@article{portnoy1986,
  title = {On the Central Limit Theorem in {$R^p$} When $p\to\infty$},
  author = {Portnoy, Stephen},
  year = 1986,
  month = nov,
  journal = {Probability Theory and Related Fields},
  volume = {73},
  number = {4},
  pages = {571--583},
  issn = {0178-8051, 1432-2064},
  doi = {10.1007/BF00324853},
  urldate = {2025-10-15},
  langid = {english}
}

@article{tang2025,
  title = {Laplace and Saddlepoint Approximations in High Dimensions},
  author = {Tang, Yanbo and Reid, Nancy},
  year = 2025,
  month = aug,
  journal = {Bernoulli},
  volume = {31},
  number = {3},
  issn = {1350-7265},
  doi = {10.3150/24-BEJ1758},
  urldate = {2025-10-15},
  langid = {english}
}

@book{vershynin2018,
  title = {High-{{Dimensional Probability}}: {{An Introduction}} with {{Applications}} in {{Data Science}}},
  shorttitle = {High-{{Dimensional Probability}}},
  author = {Vershynin, Roman},
  year = 2018,
  month = sep,
  publisher = {Cambridge University Press},
  doi = {10.1017/9781108231596},
  urldate = {2026-09-17},
  copyright = {https://www.cambridge.org/core/terms},
  isbn = {978-1-108-23159-6 978-1-108-41519-4}
}

@article{withers2026,
  title = {Edgeworth {{Expansions When}} the {{Parameter Dimension Increases}} with {{Sample Size}}},
  author = {Withers, Christopher Stroude},
  year = 2026,
  month = apr,
  journal = {Econometrics},
  volume = {14},
  number = {2},
  pages = {21},
  issn = {2225-1146},
  doi = {10.3390/econometrics14020021},
  urldate = {2026-09-22},
  langid = {english}
}

@article{zhang2012,
  title = {The {{Best Rank-1 Approximation}} of a {{Symmetric Tensor}} and {{Related Spherical Optimization Problems}}},
  author = {Zhang, Xinzhen and Ling, Chen and Qi, Liqun},
  year = 2012,
  month = jan,
  journal = {SIAM Journal on Matrix Analysis and Applications},
  volume = {33},
  number = {3},
  pages = {806--821},
  issn = {0895-4798, 1095-7162},
  doi = {10.1137/110835335},
  urldate = {2025-10-22},
  langid = {english}
}

@article{zhang2019_spa,
  title = {Fast {{Likelihood-Based Inference}} for {{Latent Count Models Using}} the {{Saddlepoint Approximation}}},
  author = {Zhang, W. and Bravington, M. V. and Fewster, R. M.},
  year = 2019,
  month = sep,
  journal = {Biometrics},
  volume = {75},
  number = {3},
  pages = {723--733},
  issn = {0006-341X, 1541-0420},
  doi = {10.1111/biom.13030},
  urldate = {2026-09-29},
  copyright = {https://academic.oup.com/journals/pages/open\_access/funder\_policies/chorus/standard\_publication\_model},
  langid = {english}
}

@article{zhilova2022,
  title = {New {{Edgeworth-type}} Expansions with Finite Sample Guarantees},
  author = {Zhilova, Mayya},
  year = 2022,
  month = oct,
  journal = {The Annals of Statistics},
  volume = {50},
  number = {5},
  issn = {0090-5364},
  doi = {10.1214/22-AOS2192},
  urldate = {2025-08-02},
  langid = {english}
}

\end{document}